\documentclass[oneside, 12pt]{amsart}
\usepackage{amscd, amssymb, amsmath, mathrsfs}
\usepackage[english]{babel}
\usepackage{booktabs}
\usepackage{tikz-cd}
\usepackage{url}
\usepackage[pdftex, colorlinks=true,  citecolor=blue, linkcolor=blue, linktocpage=true]{hyperref}

\newcommand\mylabel[1]{\label{#1}\marginpar{\vspace{-1ex}\medskip\medskip\footnotesize \tt #1}}
\renewcommand\mylabel[1]{\label{#1}}
\newcommand{\mydate}{
\number\day\space
\ifcase\month \or January\or February\or March\or April\or May\or June\or July\or August\or September\or October\or November\or December\fi 
\space\number\year}

\DeclareUrlCommand\arXiv{\urlstyle{same}}

\newtheorem{theorem}{Theorem}[section]
\newtheorem*{maintheorem}{Theorem}
\newtheorem{lemma}[theorem]{Lemma}
\newtheorem{proposition}[theorem]{Proposition}
\newtheorem{corollary}[theorem]{Corollary}

\theoremstyle{definition}
\newtheorem{definition}[theorem]{Definition}

\newtheorem*{acknowledgement}{Acknowledgement}

\theoremstyle{remark}

\newcommand{\ZZ}{\mathbb{Z}}

\newcommand{\RR}{\mathbb{R}}

\renewcommand{\AA}{\mathbb{A}}
\newcommand{\GG}{\mathbb{G}}

\newcommand{\ideala}{\mathfrak{a}}

\newcommand{\shA}{\mathscr{A}}

\newcommand{\shE}{\mathscr{E}}
\newcommand{\shF}{\mathscr{F}}
\newcommand{\shG}{\mathscr{G}}

\newcommand{\shM}{\mathscr{M}}

\newcommand{\shT}{\mathscr{T}}

\newcommand{\shV}{\mathscr{V}}
\newcommand{\shX}{\mathscr{X}}

\newcommand{\catC}{\mathcal{C}}

\newcommand{\Aff}{\text{\rm Aff}}
\newcommand{\alg}{\text{\rm alg}}

\newcommand{\Aut}{\operatorname{Aut}}

\newcommand{\Br}{\operatorname{Br}}

\newcommand{\can}{\text{\rm  can}}

\newcommand{\cl}{\operatorname{cl}}

\newcommand{\crs}{\text{\rm crs}}

\newcommand{\End}{\operatorname{End}}

\newcommand{\GL}{\operatorname{GL}}

\newcommand  {\Grass}{\operatorname{Grass}}

\newcommand{\Hom}{\operatorname{Hom}}
\newcommand{\Hilb}{\operatorname{Hilb}}

\newcommand{\id}{{\operatorname{id}}}

\newcommand{\Isom}{\operatorname{Isom}}

\newcommand{\Kernel}{\operatorname{Ker}}

\newcommand{\lra}{\longrightarrow}

\newcommand{\Mat}{\operatorname{Mat}}

\newcommand{\maxid}{\mathfrak{m}}

\newcommand{\N}{\operatorname{N}}

\newcommand{\primid}{\mathfrak{p}}
\renewcommand{\O}{\mathscr{O}}

\newcommand{\ob}{\operatorname{ob}}
\newcommand{\op}{\text{\rm op}}

\newcommand{\Pic}{\operatorname{Pic}}
\newcommand{\PGL}{\operatorname{PGL}}

\newcommand{\pr}{\operatorname{pr}}

\newcommand{\quadand}{\quad\text{and}\quad}

\newcommand{\ra}{\rightarrow}

\newcommand{\rank}{\operatorname{rank}}

\newcommand{\Set}{{\text{\rm Set}}}

\newcommand{\Sh}{{\operatorname{Sh}}}

\newcommand{\SL}{\operatorname{SL}}

\newcommand{\Spec}{\operatorname{Spec}}

\newcommand{\Sym}{\operatorname{Sym}}

\newcommand  {\uEnd}    {\underline{\operatorname{End}}} 
\newcommand  {\uAut}    {\underline{\operatorname{Aut}}}
\newcommand  {\uHom}    {\underline{\operatorname{Hom}}}

\newcommand  {\QCoh}    {\operatorname{QCoh}} 
\newcommand  {\taut}    {\text{{\rm taut}}} 
\newcommand	{\azu}{{\text{\rm azu}}}
\newcommand  {\conj}    {\operatorname{conj}} 

\newcommand{\quiver}{\Delta}
\newcommand{\pathalg}[2][\quiver]{#2[#1]}
\newcommand{\Endomor}{  E}
\newcommand{\tautsheaf}[1][V']{\shF^\taut_{#1}}
\newcommand{\tautsplit}[1][V']{\psi^\taut_{#1}}

\def\mydate{\number\day\space\ifcase\month \or January\or February\or March\or 
April\or May\or June\or July\or
August\or September\or October\or November\or December\fi \space\number\year}

\DeclareFontFamily{U}{wncy}{}
\DeclareFontShape{U}{wncy}{m}{n}{<->wncyr10}{}
\DeclareSymbolFont{mcy}{U}{wncy}{m}{n}
\DeclareMathSymbol{\SH}{\mathord}{mcy}{"58}
\begin{document}

\title[Rational Points in coarse moduli spaces]
      {Rational points in coarse moduli spaces and twisted representations}

\author[Fabian Korthauer]{Fabian Korthauer}
\address{Heinrich Heine University D\"usseldorf, Faculty of Mathematics and Natural Sciences, Mathematical Institute, 40204 D\"usseldorf, Germany}
\curraddr{}
\email{Korthauer.maths@gmx.de}

\author[Stefan Schr\"oer]{Stefan Schr\"oer}
\address{Heinrich Heine University D\"usseldorf, Faculty of Mathematics and Natural Sciences, Mathematical Institute, 40204 D\"usseldorf, Germany}
\curraddr{}
\email{schroeer@math.uni-duesseldorf.de}

\dedicatory{12 January 2025}
 
% 14D22 Fine and coarse moduli spaces
% 14D23 Stacks and moduli problems
% 14A20 Generalizations (algebraic spaces, stacks)
% 16G10 Representations of associative Artinian rings
% 16G20 Representations of quivers and partially ordered sets
% 16H05 Separable algebras (e.g., quaternion algebras, Azumaya algebras, etc.)

% Possible journals:
% Moduli. Because subject. Editor Hoskins!
% Ann. Inst. Fourier (Grenoble). Because Hoskins-Schaffhauser paper; editor ?
% Representation Theory. Editor?

\begin{abstract}
We  study  moduli spaces and moduli stacks  for representations of associative algebras
in Azumaya algebras, in rather general settings. 
We do not impose any stability condition and work over arbitrary ground rings, 
but restrict   attention  to the so-called Schur representations,
where  the only automorphisms are scalar multiplications. 
The stack comprises   twisted representations, which are representations that  live on the gerbe of
splittings for the Azumaya algebra. 
Such generalized spaces and  stacks    appear naturally: For any rational point on the classical coarse moduli space
of  matrix  representations, the machinery of non-abelian cohomology produces
a modified moduli problem for which the point acquires geometric origin. The latter are  given by 
representations in Azumaya algebras.
\end{abstract}

\maketitle
\tableofcontents

%===========================================================
\section*{Introduction}
\mylabel{introduction}

In algebraic geometry, the discrepancy between a \emph{coarse moduli space} and a  \emph{fine moduli stack}
is an endless source of puzzlement, wonder  and  research.
Loosely speaking, a ``moduli space'' is a geometric object that parameterizes a given class of other geometric objects,
in an ``optimal way''.
For example, the Hilbert scheme $\Hilb_{X/k}$ parameterizes the closed subschemes $Z$ of a given projective scheme $X$, say
over a ground field $k$.
This is the prime example of a \emph{fine moduli space}, which   represents the corresponding functor of families of closed subschemes.
If $X$ is merely proper, not necessarily having an ample invertible sheaf, $\Hilb_{X/k}$ still exists as an \emph{algebraic space},
which are certain sheaves $(\Aff/k)\ra (\Set)$   on the category of affine schemes over $k$  that are   closely related to schemes. 
Despite this additional layer of technicalities, $\Hilb_{X/k}$ remains a fine moduli space.

In most other cases of interest, the situation is   more difficult, due to non-trivial automorphisms of objects. 
Here the prime example is the coarse moduli space   $M_g$ that  parameterizes smooth curves of genus $g\geq 2$.
Although it is true that $M_g(\Omega)$ corresponds to isomorphism classes of such curves defined over algebraically closed fields $\Omega$, this does not
hold true for general  fields $k$, let alone    rings $R$. Similar things can be said about  the
Picard scheme $\Pic_{X/k}$. In fact, in presence of non-trivial automorphisms a fine moduli space
\emph{never  exists}, a fact that cannot be over-emphasized.

To cope with such issues, one has to replace schemes $(X,\O_X)$ or  algebraic spaces  $X:(\Aff/k)\ra (\Set)$ by stacks $\shX\ra (\Aff/k)$, which 
are certain fibered groupoids.
This turns all our players into categories. In particular, a scheme or an algebraic space $X$ 
becomes the  stack $\shX=(\Aff/X)$ via the comma construction. A fundamental insight of Deligne and Mumford
\cite{Deligne; Mumford 1969} and  Artin \cite{Artin 1974}  was that the geometric intuition stemming from
schemes still works very well in the realm of \emph{Deligne--Mumford stacks} or even  \emph{Artin stacks}.
For comprehensive general presentations, see \cite{Laumon; Moret-Bailly 2000}, \cite{Olsson 2016} and \cite{SP}.

Our initial motivation was to generalize 
a beautiful work of Hoskins and Schaff\-hauser \cite{Hoskins; Schaffhauser 2020} on rational points on moduli of quiver representations 
into an abstract setting.
Recall that a \emph{quiver} $\quiver$ is a finite directed graph, with loops and multiple arrows allowed.
A \emph{quiver representation} attaches to each vertex $i\in \quiver$ a finite-dimensional vector space $E_i$, 
and to each arrow $i\stackrel{\alpha}{\ra} j$ in $\quiver$
a linear map $f_\alpha:E_i\ra E_j$.
Let $M_{\quiver/k}$ be the resulting coarse moduli space  of stable quiver representations  over some ground field $k$,
with fixed dimension vector and  with respect  to some stability condition, both omitted from notation.
The construction involves the formation of a GIT quotient \cite{Mumford; Fogarty; Kirwan 1993}, 
and consequently the rational points $x\in M_{\quiver/k}(k)$ do not necessarily come
from an actual quiver representation $(E_i\mid f_\alpha)_{i,\alpha\in\quiver}$. 
In \cite{Hoskins; Schaffhauser 2020} the obstructions are analyzed in an explicit, down-to-earth way, 
and a geometric interpretation in terms of modified data are given, of course involving the \emph{Brauer group} $\Br(k)$.

We put ourselves into a general, abstract and categorified situation, which encompasses all sorts of applications:
Let us work in a topos $\shT$ with final object $S$, and a central extension of group objects $1\ra N\ra G\ra H\ra 1$,
with $H$ acting freely on some  object $X$, with ensuing quotient $Q=X/H^\op$.
Intuitively, we think of $X$ as a space that ```over-parameterizes'' a certain class of geometric objects.
The $H$-action corresponds to the isomorphism relation, whereas $N$ reflects   automorphisms.
In turn, one may say that $Q$ ``under-parameterizes'' our class of geometric objects.

We  now say that an $S$-valued point $g:S\ra Q$ is \emph{of geometric origin} if it admits a lifting $\tilde{g}:S\ra X$.
One obtains an obstruction map 
$\ob:Q(S)\ra H^2(S,N)$. By non-abelian cohomology and twisted forms we  construct, 
with fixed $g\in Q(S)$ and in a canonical way, a ``modified moduli problem'' 
$\tilde{X}$ and   $1\ra N\ra \tilde{G}\ra \tilde{H}\ra 1$, without changing the quotient $Q$.  Our first  result is:

\begin{maintheorem}
(See Thm.\ \ref{twist geometric origin})
With respect to the quotient  $Q=\tilde{X}/\tilde{H}^\op$ 
for the modified moduli problem, the given $S$-valued point $g\in Q(S)$   acquires geometric origin.
\end{maintheorem}

In one form or another, this appears  to be well-known (compare  
 the monograph of Skorobogatov \cite{Skorobogatov 2001}, Section 2.2).

We     apply our theory  to \emph{representations of   associative algebras} $\Lambda$ over   ground rings $R$.
Note that this encompasses representations of  groups $\Gamma$  via the \emph{group algebra} $\Lambda=R[\Gamma]$,
representations of quivers $\quiver$ via the \emph{path algebra} $\Lambda=\pathalg{R}$. 
This apparently   applies to     coherent sheaves  
on projective schemes   as well, by an ingenious construction of \'Alvarez-C\'onsul and  King \cite{Alvarez-Consul; King 2007}
relying on  \emph{Kronecker modules}; it would be very interesting to pursue this further. 

For simplicity we assume that the associative algebra $\Lambda$ is finitely presented.
To work with $H=\PGL_n$ and $G=\GL_n$, we have to restrict to representations
whose endomorphism ring comprises only scalars. Note that there is no consistent designation for this  important property.
Following Derksen and Weyman (\cite{Derksen; Weyman 2011}, Section 2.3) we like to call them \emph{Schur representations}.
One also finds the terminology  ``brick'' (\cite{Assem; Simson; Skowronski 2006}, Chapter VII, Definition 2.4) and
``stably indecomposable''  (\cite{Kraft; Riedtmann 1986}, Section 2.6).
In all other respects we work in far greater generality than customary in geometric invariant theory, 
\emph{without any stability conditions and over arbitrary ground rings}.
Note that over fields, the Schur representations contain all geometrically stable representations.

According to Grothendieck's fundamental insight \cite{GB}, the twisted forms $\tilde{H}$ of the group scheme $H=\PGL_n$  correspond
to \emph{Azumaya algebras}. Consequently, it becomes imperative  to  consider representation of $\Lambda$ not only in matrix algebras $\Mat_n(R)$,
but in general Azumaya algebras  $\Lambda^\azu$ of degree $n\geq 1$. Our second main result is:

\begin{maintheorem}
(See Thm.\ \ref{space  schur representations}) Suppose $\Lambda$ is finitely presented as an  associative $R$-algebra,
and $\Lambda^\azu$ is an Azumaya algebra of degree $n\geq 1$. Then the functor $X=X^{\Lambda^\azu}_{\Lambda/R}$ of Schur representations
is   representable by a quasiaffine scheme of finite presentation. Moreover, the group scheme $H=\Aut_{\Lambda^\azu/R}$ acts freely,
and the quotient $Q=X/H^\op$ is an algebraic space that is of finite presentation.
\end{maintheorem}

For the path algebra $\Lambda=\pathalg{R}$ of a quiver we also give an independent argument using Grassmann varieties and vector bundles.
Note that \emph{algebraic spaces}  are   generalizations of schemes introduced by Artin 
(\cite{Artin 1971},  \cite{Knutson 1971}, \cite{Laumon; Moret-Bailly 2000}, \cite{Olsson 2016}) that allow
the formation of quotients which usually do no exist as schemes. Moreover, they are indispensable
in the very definition of Artin stacks. Note that algebraic spaces like $Q=X/H^\op$ are frequently non-separated, and
easily become non-schematic  (compare \cite{Schroeer 2022}).
 
Our third main result is a description of the   moduli stack $[X/G/Q]$, where $G=U_{\Lambda^\azu/R}$ is the 
group scheme of units  of the Azumaya algebra $\Lambda^\azu$. Quotient stacks are constructed in the geometric
language of principal homogeneous spaces. Our description is in a representation-theoretic way, much closer to the moduli problem at hand:
We define a \emph{twisted Schur representation} of $\Lambda$ in $\Lambda^\azu$ over an affine scheme $V$
as a pair $(\shE',\rho')$, where $\shE'$ is a locally free sheaf of rank $n=\deg(\Lambda^\azu)$ and weight $w=1$
on the  gerbe $V'$ of splittings for $\Lambda^\azu\otimes_R\O_{V'}$, and $\rho':\Lambda\otimes_R\O_{V'}\ra \uEnd(\shE)$
is a Schur representation.
This of course relies on the notion of \emph{twisted sheaves}, 
which where   introduced by C\u{a}ld\u{a}raru  (\cite{Caldararu 2000}, \cite{Caldararu 2002}),
de Jong \cite{de Jong 2006} 
and Lieblich (\cite{Lieblich 2004} and  \cite{Lieblich 2007}), and have attracted tremendous interest in the past decades. 
Earlier, Edidin, Hassett, Kresch and  Vistoli \cite {Edidin; Hassett; Kresch; Vistoli 1999}   already established  a deep link between   
the existence of Azumaya algebras representing a given Brauer class  and internal properties of  the  Artin stacks of splittings. 
Our third main result is:

\begin{maintheorem}
(See Thm.\ \ref{comparison functor equivalence}) The stack of twisted Schur representations
$$
\shM_{\Lambda/R}^{\Lambda^\azu}=\{(V,\shE',\rho')\}
$$
is equivalent to the quotient stack $[X/G^\op/Q]$, where $X=X_{\Lambda/R}^{\Lambda^\azu}$ is the quasiaffine scheme that represents the functor
of Schur representations of $\Lambda$ in $\Lambda^\azu$.
\end{maintheorem}

In other words,  $\shM_{\Lambda/R}^{\Lambda^\azu}$ is the ``true'' moduli stack for representations of associative rings $\Lambda$
in Azumaya algebras $\Lambda^\azu$. The algebraic space $Q=X/H^\op$ is the coarse moduli space over $S=\Spec(R)$.
Given some $g\in Q(S)$, we now can   determine when it has geometric origin,
and  describe the modified moduli problem $Q=\tilde{X}/\tilde{H}$, for which $g$ has acquired geometric origin:

\begin{maintheorem}
(See Thm.\ \ref{characterization geometric origin} and \ref{representation-theoretic content}) 
If the non-abelian cohomology set $H^1(S,\GL_n)$ is a singleton
and the category $\shM_{\Lambda/R}^{\Lambda^\azu}(S)$ is non-empty, then $g$ has geometric origin for $Q=X/H^\op$.
In any case, the twisted form $\tilde{X}$ is the scheme of Schur representations of $\Lambda$
in the Azumaya algebra $\tilde{\Lambda}^\azu$ obtained by twisting $\Lambda^\azu$. 
\end{maintheorem}

Our fourth main  result deals with the   \emph{tautological sheaf} $\shT_\shM$ on the moduli stack
$\shM=\shM_{\Lambda/R}^{\Lambda^\azu}$. Locally, this is a Hom sheaf between two locally free sheaves
on the gerbe of splittings of weight one.  Its endomorphism algebra has weight zero, and thus yields
an Azumaya algebra $\shA_Q$ over the algebraic space   $Q=X/H^\op$, 
which contains a smaller Azumaya algebra $\shA^0_Q$ as the commutant of $\Lambda^\azu$.
Their classes in the Brauer group $\Br(Q) \subset H^2(Q,\GG_m)$ have the following meaning:

\begin{maintheorem}
(See Thm.\ \ref{tautological classes}) 
In the cohomology group $H^2(Q,\GG_m)$ of the algebraic space $Q=X/H^\op$, we have
$[\shA_Q]=\partial[X]$ and $[\shA^0_Q]=[\Lambda^\azu\otimes_R\O_Q]$.
\end{maintheorem}

Here $[X]\in H^1(Q,H)$ is the torsor class for the quotient map $X\ra Q$, and $\partial[X]$ is its Brauer class stemming
from the non-abelian coboundary map.

\medskip
The paper is structured as follows:
In Section \ref{Points geometric origin} we collect some generalities on torsors, introduce the points of geometric origin, and
the technique of modifying moduli problems. Section \ref{Quotient stacks} contains generalities on quotient stacks and gerbes.
In Section \ref{Spaces representations} we introduce the notion of Schur representations of an associative algebra $\Lambda$
in another associative algebra $\Lambda'$ over arbitrary ground rings $R$, and establish   representability results
for the functor $X=X^{\Lambda^\azu}_{\Lambda/R}$ and quotient $Q=X/H^\op$.
In Section  \ref{Quiver representations} we give an  alternative approach for quiver representations.
Section \ref{Gerbe splittings} contains a detailed and explicit discussion of   splitting gerbes for Azumaya algebras.
The central part of the paper is Section \ref{Stacks schur},
where we introduce the Artin stack $\shM=\shM^{\Lambda^\azu}_{\Lambda/R}$ of twisted Schur representations,
and show that it is equivalent to a quotient stack $[X/G/Q]$. In Section \ref{Modifying moduli}, 
we apply  our findings to understand the modification
of the moduli problem of Schur representations.
In Section \ref{Tautological sheaves} we introduce the tautological sheaf $\shT_\shM$ on the moduli stack $\shM$, 
and use it to explain various Azumaya algebras
and Brauer classes.

\begin{acknowledgement}
This research was also conducted in the framework of the   research training group
\emph{GRK 2240: Algebro-Geometric Methods in Algebra, Arithmetic and Topology},  which is  funded
by the Deutsche Forschungsgemeinschaft. The first author received support  by the GRK 2240 via a  
start-up funding grant. 

\end{acknowledgement}

%===========================================================
\section{Points of geometric origin}
\mylabel{Points geometric origin}

In this section we review the theory of torsors and twisting, and introduce the notation of $S$-valued
points of geometric origin.
We will freely use the language of sites, topoi, stacks, torsors, and  gerbes, 
but sometimes recall and comment on  some relevant issues along the way.
For   details we refer to  
\cite{Grothendieck 1955}, \cite{SGA 4a}, \cite{Giraud 1971},  \cite{Laumon; Moret-Bailly 2000}, \cite{Olsson 2016},
\cite{Skorobogatov 2001}.
 
Let $\catC$ be a site   and   $\shT=\Sh(\catC)$ be the ensuing topos of sheaves.
To simplify exposition, we assume that   $\catC$ contains a terminal object $S$, all fiber products exist, and 
all representable functors $\catC\ra(\Set)$ satisfy the sheaf axiom.
The latter  gives, via the Yoneda Lemma,  a fully faithful embedding $\catC\subset\shT$, and in particular we can regard $S$ as   terminal object in $\shT$.
Throughout,   all objects belong to this $\shT$ if not said otherwise.

Let $G$ be a group object, acting on another object $Z$. For each $G$-torsor $P$, we thus obtain a \emph{twisted form}
$$
{}^P\!Z= P\wedge^GZ = G\backslash (P\times Z) = (P\times Z)/G^\op 
$$
of $Z$. Here $G$ acts diagonally on the product $P\times Z$, and $G^\op$ denotes the opposite group, which indeed acts  from the right rather than the left.
Note that the construction is functorial in $Z$, and commutes with products.
In particular, if $Z$ is endowed with some algebraic structure respected by the $G$-action,
the twisted form ${}^PZ$ inherits the algebraic structure.
In the special case $Z=G$, on which $G$  acts via  conjugation, we thus get a new  group object ${}^P\!G$. 
More generally, if  $Z$ is endowed with an object of  operators $\Omega$ 
(in the sense of  \cite{A 1-3}, Chapter I, \S3, No.\ 1), and $G$ acts on $Z$ and $\Omega$ in a compatible way,
that is, $\Omega\times Z\ra Z$ is $G$-equivariant, then ${}^P\!Z$ is endowed with  operators ${}^P\Omega$.

Suppose now that $Z=P$, and write $P_0$ for the underlying object  where the $G$-action is omitted, 
giving  an inclusion $\Aut_{P/S}\subset\Aut_{P_0/S}$.
Combining \cite{Hilario; Schroeer 2023}, Section 9
and \cite{Schroeer; Tziolas 2023}, Lemma 3.1, we have a canonical homomorphism
$$
{}^P\!G\times {}^P\!G^\op\lra {}^P\!\Aut_{G_0/S}=\Aut_{P_0/S}
$$
stemming from the action of $G\times G^\op$ on $G$   by left-right multiplication $(g_1,g_2)\cdot g= g_1gg_2$.
This endows the underlying object $P_0$ of the $G$-torsor $P$ with the additional structures
of a ${}^P\!G$-torsor and a ${}^P\!G^\op$ torsor. In fact,  ${}^P\!G^\op$ becomes the automorphism group object for the $G$-torsor $P$,
in other words
\begin{equation}
\label{automorphisms torsor}
\Aut_{P/S}= {}^P\!G^\op.
\end{equation}
Write $(\text{$G$-Tors}/S)$ and $(\text{${}^P\!G$-Tors}/S)$ for the groupoids of torsors for $G$ and ${}^P\!G$, respectively.
We then  get the \emph{torsor translation functor}
$$
(\text{${}^P\!G$-Tors}/S)\lra (\text{$G$-Tors}/S),\quad T\longmapsto P\wedge^{{}^P\!G}T,
$$
where the $G$-action is via $P$. This is indeed well-defined, 
because  in light of \eqref{automorphisms torsor} the left $G$-action on $P$ commutes
with the right ${}^P\!G$-action. It  sends the trivial torsor $T={}^P\!G$ to the torsor $P$, which   usually is non-trivial.
Likewise, we have a functor in the other direction
$$
(\text{$G$-Tors}/S)\lra (\text{${}^P\!G$-Tors}/S),\quad T\longmapsto  \Isom_G(P,T),
$$
where the ${}^P\!G$-action is via $P$, again given by  \eqref{automorphisms torsor}.
This sends $T=P$ to the trivial torsor $\Isom_G(P,P)=G\cdot\id_P$.
The above functors are mutually inverse equivalences, where the adjunction maps  
$$
T\lra\Isom_G(P,P\wedge^{{}^P\!G}T)\quadand   P\wedge^{{}^P\!G} \Isom_G(P,T)\lra T
$$
are given by $t\mapsto (p\mapsto p\wedge t)$ and $p\wedge f\mapsto f(p)$, respectively.
On the sets of isomorphism classes, we get mutually inverse \emph{torsor translation maps}
\begin{equation}
\label{torsor tranlation map}
H^1(S,{}^P\!G)\lra H^1(S,G)\quadand H^1(S,G)\lra H^1(S,{}^P\!G)
\end{equation}
of sets. Note that in general these maps    do not respect the distinguished points.
 
Suppose now that $G\ra H$ is a homomorphism of group objects. For each $G$-torsor $P$ we get an induced $H$-torsor 
$H\wedge^GP$.
If $G\ra H$ is an epimorphism with kernel $N$, this becomes $N\backslash P=P/N^\op$.
Conversely, for general $G\ra H$, we say that an $H$-torsor $T$ admits a \emph{reduction of structure group} if there is some $G$-torsor $P$
with  $H\wedge^GP\simeq T$. If the homomorphism is a monomorphism $G\subset H$, this simply means that $T/G^\op$ admits a section.
Reduction  to $G=\{e\}$ boils down to the triviality of the $H$-torsor.   

Now let $X$ be an object, and $H$ be a group object acting freely on $X$, and write 
$$
Q=H\backslash X=X/H^\op
$$
for the resulting quotient.  Note that the canonical map $X(S)\ra Q(S)$ is not surjective in general, because the formation of quotients involves
sheafification. 
The objects $H_Q=H\times Q$ and $X$ are endowed with canonical morphism to $Q$, via the second projection and the quotient map, respectively.
In turn, we may view them as objects in the \emph{comma category} $\shT_{/Q}$, comprising pairs $(U,g)$ with $U$
an object and $g:U\ra Q$ a morphism from $\shT$.
This is again a   topos (\cite{SGA 4a}, \'Expose IV, Theorem 1.2).
Furthermore, $H_Q$ has the structure of  a group object in $\shT_{/Q}$, with $X$ as an  $H_Q$-torsor.
%Reference?

Given an $S$-valued point $g:S\ra Q$ of the quotient $Q=X/H^\op$, we get an induced torsor
$g^*(X) = X\times_QS$
with respect to $H=g^*(H_Q)=(H\times Q)\times_QS$. The following locution will be useful throughout:

\begin{definition}
\mylabel{geometric origin}
In the above setting, we say that  the $S$-valued point $g:S\ra Q$ of the quotient $Q=X/H^\op$ is of \emph{geometric origin} 
if the induced $H$-torsor $g^*(X)$ is trivial.
\end{definition}

The condition simply means that the projection  $g^*(X)\ra S$ admits a section. By the cartesian diagram
$$
\begin{CD}
g^*(X)	@>>>	X\\
@VVV		@VVV\\
S	@>>g>	Q,
\end{CD}
$$
this also means that $g\in Q(S)$ lies in the image of $X(S)$. This of course is automatic if $X\ra Q$ admits a section,
but in general the map $X(S)\ra Q(S)$ is far from surjective.

A first instructive example arises from the   polynomial ring $R=\RR[u]$ 
and the   extension $A=R[t]/(f)$ given by the polynomial $f(t)=t^2 +ut+1$. Write $X$ and $Q$ for the respective spectra of the localizations
$A_\Delta$ and $R_\Delta$ with respect to the discriminant $\Delta=u^2-4$.
Then $X\ra Q$ is a finite \'etale Galois covering with  Galois group $H=\ZZ/2\ZZ$, and we see that 
an  $\RR$-valued point $s\in Q$ given by $u=\lambda$
is   of  geometric origin if and only if $|\lambda|>2$.

Our motivation for  the terminology stems from  the following:
Think  of $X$ as a space that  over-parameterizes   \emph{geometric objects} stemming from some \emph{moduli problem}, 
in such a way  that the $H$-orbits 
correspond to the \emph{isomorphism classes}.
One may regard $Q$ as a \emph{coarse moduli space},   the $S$-valued points $g:S\ra Q$ as    $S$-families of isomorphism classes,
and  the $S$-valued points $\tilde{g}:S\ra X$ as  the more significant $S$-families of   geometric objects. 
The map $X(S)\ra Q(S)$ is not necessarily surjective,
and coping with this failure is the main topic of this paper.  
In praxis, $H$ often sits in an exact sequence $1\ra N\ra G\ra H\ra 1$, where $N$ reflects the \emph{automorphisms} in the moduli problem,
a situation studied  in the next section. 

The following example elucidates what we have in mind:
Fix a ground ring $R$, endow  the category $\catC=(\Aff/R)$ of affine schemes   with the fppf topology, and 
consider the functor $Y(A)=\Mat_n(A)\times\Mat_n(A)$ of matrix pairs. The projective linear group $H=\PGL_n$ acts via simultaneous conjugation,
and we consider the subfunctor $X\subset Y$ of matrix pairs whose stabilizer is trivial.
Note that $G=\GL_n$ and $N=\GG_m$ give the short exact sequence  $1\ra N\ra G\ra H\ra 1$,
that the functor  $X$ is a quasiaffine scheme  (see Theorem \ref{space schur representations} below),   
and that the quotient $X/H^\op$ can be seen as the coarse moduli space
of    certain $n$-dimensional  representation of the associative algebra $\Lambda_0=R\langle t_1,t_2\rangle$, namely 
those  with only scalars as automorphisms.
Also note that if  $R=K$ is a field, $\Lambda_0$ is the proto-typical example of a \emph{wild algebra}. In other words,   its
representation theory ``contains'' the representation theory  
of every   finite-dimensional $K$-algebra $\Lambda$ (compare   \cite{Kirillov 2016}, Section 7.1).

Back to the general setting. 
We now fix an $S$-valued point $g:S\ra Q$, and consider the resulting $H$-torsor $P=g^*(X)=X\times_QS$. With respect to the conjugation
 action of $H$ on itself, 
we obtain a twisted form $\tilde{H}=P\wedge^HH$.
Recall that $X$ is a torsor with respect to $H_Q$, and the same holds for the pull-back $P_Q=P\times Q$. 
We now apply the above considerations in the comma category $\shT_{/Q}$,
and regard the twisted form 
$$
\tilde{X} = \Isom_{H_Q}( P_Q,X)
$$

as a torsor for $\tilde{H}_Q=({}^P\!H)_Q={}^{P_Q}\!(H_Q)$. By construction, both quotients $X/H^\op$ and $\tilde{X}/\tilde{H}^\op$ coincide with
the terminal object  $Q$ from the comma category, and thus  
$$
X/H^\op=Q=\tilde{X}/\tilde{H}^\op.
$$
Our first main result is the following  foundational fact, now almost a triviality:

\begin{theorem}
\mylabel{twist geometric origin}
With respect to the quotient  $Q=\tilde{X}/\tilde{H}^\op$, the given $S$-valued point $g:S\ra Q$ has acquired geometric origin.
\end{theorem}

\proof
Making a base-change, it  suffices to treat the case that $S=Q$, $g=\id_Q$ and $P=X$.
Then $\tilde{X}$ is just the image of $T=P$ under the torsor translation map 
$$
H^1(S,H)\lra H^1(S,{}^P\!H),\quad T\longmapsto \Isom_H(P,T).
$$
We already remarked that the object $\Isom_H(P,P)$ admits the identity as section.
\qed

\medskip
Intuitively, one should regard the passage from $X$ to  $\tilde{X}$ as  \emph{modifying} the   moduli problem  \emph{without} changing the coarse moduli space.
Thus $g:S\ra Q$, which   is just an $S$-family of isomorphism classes for the initial moduli problem,   
is promoted  to an $S$-family of geometric objects for the modified moduli problem.
\emph{For a given concrete moduli problem $X$, this raises the question to 
understand the geometric meaning of the new moduli problem $\tilde{X}$.} 
In the following sections we will obtain   answers for certain moduli of representations.

Also note that the set $Q(S)$  comes with a partition
$$
Q(S) = \dot{\bigcup_{[P]\in H^1(S,H)}} Q(S)_P
$$
where the  $Q(S)_P$   comprise the $g:S\ra Q$ for which $g^*(X)$ is isomorphic to $P$. 
Note that these are exactly the $S$-valued points which acquire geometric origin with respect to the quotient $Q=\tilde{X}/\tilde{H}^\op$, 
where $\tilde{X}$ and $\tilde{H}$ are obtained by twisting $X$ and $H$ with $P$.
In the context of quiver moduli, 
the above partition was already obtained by Hoskins and Schaffhauser, where  the quotient
map $Q(S)\ra H^1(S,H)$ was called the \emph{type map}  (\cite{Hoskins; Schaffhauser 2020}, Theorem 1.1 and Remark 3.5).

%===========================================================
\section{Quotient stacks}
\mylabel{Quotient stacks}

\newcommand{\A}{N}
We keep the set-up of the preceding section.
Let $X$ be an object from the topos $\shT=\Sh(\catC)$, and $H$ be a group object, acting freely on $X$ with quotient $Q=X/H^\op$.
In this section, we additionally assume that $H$ sits in some short exact sequence
\begin{equation}
\label{extension}
1\lra \A\lra G\lra H\lra 1
\end{equation}
of group objects, where we assume that   the epimorphism   $G\ra H$ locally admits sections, not necessarily respecting  the group laws.
Note that this holds with $G=\GL_n$ and $H=\PGL_n$,  but fails for $G=H=\GG_m$ and taking powers by $d\geq 2$.
We now ask whether the $H$-torsor $g^*(X)=X\times_QS$, for fixed  $g\in Q(S)$,  admits a   reduction 
of structure group with respect to the epimorphism $G\ra H$. This of course holds if the  torsor is trivial, but is in general a much weaker condition.

If $\A$ is non-trivial, the $G$-action   on $X$ is non-free. Instead of forming the quotient $X/G^\op=X/H^\op=Q$, which loses
too much information, we consider  the category
$$
[X/G^\op/Q]=\{(U,g,P,f)\mid f:P\ra X_U\} 
$$
whose  objects are quadruples where  $U$ is an object from $\catC$, and $g:U\ra Q$ is a morphism in $\shT$, and $P$ is a $G_U$-torsor,
and $f:P\ra X_U$ is a morphism that is  $G_U$-equivariant. Note that the latter condition holds if and only if the
composite map $P\stackrel{f}{\ra} X_U\stackrel{\pr_1}{\ra} X$ is $G$-equivariant. Also note that $g$ is determined as the map induced
from $\pr_1\circ f$.

Morphisms $(U',g',P',f')\ra(U,g,P,f)$ are defined as  pairs $(u,p)$ of morphisms $u:U'\ra U$ and $p:P'\ra P$  
such that the diagram
$$
\begin{CD}
X_{U'}		@<f'<<	P'	@>>>	U'	@>g'>>	Q\\
@V\id\times uVV		@VVpV		@VVuV		@VV\id V\\
X_U		@<<f<	P	@>>>	U	@>>g>	Q
\end{CD}
$$
is commutative, and the resulting   $P'\ra P\times_UU'$ is $G_{U'}$-equivariant.
The category  comes with a forgetful functor
$$
[X/G^\op/Q]\lra \catC_{/Q},\quad (U,g,P,f)\longmapsto (U,g)
$$
to the comma category of objects $U\in \catC$ endowed with a structure morphism to $Q\in\shT$.

\begin{proposition}
\mylabel{quotient fibered category}
The above functor gives $[X/G^\op/Q]$ the structure of a fibered groupoid over $\catC_{/Q}$.
\end{proposition}

\proof
For each $(U,g)\in\catC_{/Q}$, a morphism $(u,p):(U,g,P',f')\ra (U,g,P,f)$ in the fiber category
has $u=\id_U$, and is thus determined by the morphism $p:P'\ra P$ of $G_U$-torsors. Any morphism of torsors is necessarily
an isomorphism, so the fiber categories are groupoids.

It remains to check that every morphism $(u,p):(U',g',P',f')\ra(U,g,P,f)$ is cartesian.
The induced morphism in $\catC_{/Q}$ is $u:U'\ra U$, and $p$ is already determined by the isomorphism
$P'\ra P\times_UU'$. The universal property of fiber products in $\shT$ reveals that $(u,p)$ is cartesian.
\qed

\medskip
We now regard the comma category $\catC_{/Q}$ as a site: A family  $((U_\lambda,g_\lambda)\ra (U,g))$
is a \emph{covering} if and only if  $(U_\lambda \ra U)$ is a covering family  for the site $\catC$.

\begin{proposition}
\mylabel{quotient gerbe}
The fibered groupoid  $[X/G^\op/Q]\ra \catC_{/Q}$ satisfies the stack axioms, and is actually a gerbe.
\end{proposition}

\proof 
Using that for $\Sh(\catC)$ all  descent data are effective, one easily shows that
every descent datum for $[X/G^\op/Q]\ra \catC_{/Q}$ is effective.
Now let $(U,g,P,f)$ be an  object over some $(U,g)$. For each  morphism $(U',g')\ra (U,g)$ in the
comma category, we choose  pullbacks  $(U',g',P',f')$ in $[X/G^\op/Q]$. Using the sheaf property of $P\in\shT$, one easily checks that
the resulting presheaf
$$
(U',g')\longmapsto \Aut \left((U',g',P',f')/(U',g')\right),
$$
satisfies the sheaf axiom. Note that the automorphism group is formed in the fiber category.
Summing up,  $[X/G^\op/Q]\ra \catC_{/Q}$ is a stack.

% {\Fvorzwei Next we check that locally the fibres of $[X/G^\op/Q]\ra \catC_{/Q}$ are non-empty. 
% If $(U,g)$ is an object of $\catC_{/Q}$, we may choose a covering family $(U_\lambda\ra U)$ 
% such that $X_{U_\lambda}\ra U_\lambda$ admits a section for each $\lambda$. 
% These sections extend to equivariant maps $f_\lambda: G_{U_\lambda}\to X_{U_\lambda}$, 
% so $(U_\lambda,g\vert U_\lambda, G_{U_\lambda}, f_\lambda)$ are objects in the fibres of $(U_\lambda,g\vert U_\lambda)$.}

It remains to check that all objects in the stack over a fixed $(U,g)$ are locally isomorphic.
Fix two  objects  $(U,g,P_1,f_1)$ and $(U,g,P_2,f_2)$, and choose a covering family $(U_\lambda\ra U)$ 
that trivializes the $G_U$-torsors $P_1$ and $P_2$. This reduces us to the special case $P_i=G_U$.
Write $f_1(e_G)= \sigma \cdot f_2(e_G)$ for some section $\sigma:U\ra H_U$.
Choose another covering family  $(U_\mu\ra U)$ such that $\sigma|U_\mu$ admits a lift to $G_{U_\mu}$. Note that here we use
our standing assumption that the epimorphism    $G\ra H$ locally admits sections, not necessarily respecting the group laws.
This reduces our problem to the case that $\sigma:U\ra H_U$ arises from some $\tilde{\sigma}:U\ra G_U$.
Then  $(\id_U,\tilde{\sigma})$ defines the desired isomorphism between $(U,g,P_1,f_1)$ and $(U,g,P_2,f_2)$.
\qed

\medskip
For each $(U,g)\in\catC_{/Q}$, we have a canonical identification of comma categories
$(\catC_{/Q})_{/(U,g)}=\catC_{/U}$. We now seek to compute the sheaf of groups 
$$
\uAut_{(U,g,P,f)/(U,g)} \subset\uAut_{P/U}={}^P\!G_U^\op,
$$
where the last identification comes from \eqref{automorphisms torsor}.
On the other hand, twisting the exact sequence \eqref{extension} yields another exact sequence $1\ra {}^P\!\A_U\ra {}^P\!G_U\ra{}^P\!H_U\ra 1$, which gives 
an inclusion 
${}^P\!\A_U^\op\subset {}^P\!G_U^\op=\uAut_{P/U}$.

\begin{proposition}
\mylabel{automorphism sheaves}
We have $\uAut_{(U,g,P,f)/(U,g)}={}^P\!\A_U^\op$ as subsheaves inside $\uAut_{P/U}$.
\end{proposition}
 
\proof
The problem is local in $U$, so it suffices to treat the case that $P=G_U$, and thus ${}^P\!\A_U^\op=\A_U^\op$.
Given a section $a:U\ra \A_U$, we write $p:G_U\ra G_U$ for the right-multiplication with $a$.
Then $(\id_U,p)$ is the desired automorphism of $(U,g,P,f)$. Conversely, each automorphism $(\id_U,p)$
yields a commutative diagram
$$
\begin{CD}
X_U	@<<<	G_U\\
@V\id VV		@VVpV\\
X_U	@<<<	G_U.
\end{CD}
$$
The map on the right takes the form $p(x)=xg$ for some $g\in G(U)$. Since the $H$-action on $X$ is free,
the above diagram ensures $g\in N(U)$.
\qed

\medskip
Let us say that $N$ is \emph{central} if the subgroup object $N\subset G$ is contained in the center $Z(G)$.
In turn, $N$ is abelian, the  conjugation action of $G$ on $N$ is trivial, and \eqref{extension} is a central extension.
In particular we have identifications
${}^P\!\A_U^\op = \A_U^\op=\A_U$. Thus for each object $(U,g,P,f)\in [X/G^\op/Q]$ we get an isomorphism
$$
\varphi_{(U,g,P,f)}:\A_U\lra \uAut_{(U,g,P,f)/(U,g)},
$$
which are compatible with respect to morphisms. In other words: 

\begin{corollary}
\mylabel{abelian gerbe}
If $N$ is central, the datum of the above maps endows the stack $[X/G^\op/Q]$ over $\catC_{/Q}$ with the structure of an $\A_Q$-gerbe.
\end{corollary}
 
Suppose   $N$ is central. For each $g:U\ra Q$,  the resulting $\N_U$-gerbe
yields a  \emph{cohomology class}
$$
\cl([g^*X/G^\op/U])=g^*\cl([X/G^\op/Q])\in H^2(U,\N_U),
$$
which are compatible with respect to pullbacks. We refer to \cite{Giraud 1971}, Chapter IV, Section 3.4 for more details.
In fact, by the geometric interpretation of low-degree cohomology,
one may view  $H^1(U,N_U)$ as the set of isomorphism classes of $N_U$-torsors,
and $H^2(U,N_U)$ as the set of equivalence classes of $N_U$-gerbes. Together with $H^0(U,N_U)=\Gamma(U,N_U)$,
these indeed form a  delta functor in degree $i\leq 2$, in the sense of \cite{Grothendieck 1957}.

For $U=S$, we obtain the so-called \emph{obstruction  map}
$$
\ob:Q(S)\lra H^2(S,\A),\quad g\longmapsto g^*\cl([X/G^\op/Q]).
$$
Under an additional assumption, this is \emph{precisely} the obstruction against having geometric origin:

\begin{theorem}
\mylabel{obstruction geometric origin}
Suppose $\A$  is central, and that the canonical map $H^1(S,N)\ra H^1(S,G)$ is surjective.  
Then $g\in Q(S)$ is of geometric origin if and only if the obstruction $\ob(g)\in H^2(S,\A)$ is zero.
\end{theorem}
 
\proof
Saying that $g:S\ra Q$ has geometric origin means that the $H$-torsor $g^*(X)$ is trivial.
Then  it obviously admits a reduction of structure group with respect to $G\ra H$, and thus $\ob(g)$ vanishes.
Conversely, suppose that $\ob(g)=0$. Then $g^*(X)$ arises from a $G$-torsor $P$. The latter comes from
an  $N$-torsor $T$, by assumption. Thus $g^*(X)$ is induced from an  $N$-torsor $T$ with respect to the composite map $N\ra G\ra H$,
which is trivial, and thus $g^*(X)$ is trivial.
\qed

\medskip
Let us now specialize the above result to the case that $G=\GL_n$ and $N=\GG_m$, and $\catC=(\Aff/R)$ is the category 
of affine schemes over a ground ring $R$, endowed with the fppf topology. The final object then is $S=\Spec(R)$.

\begin{corollary}
\mylabel{local or dedekind}
In the above setting,  suppose that the ground ring $R$ is 
\begin{enumerate}
\item semilocal,
\item or a polynomial ring in finitely many indeterminates, either over a field   or a principal ideal domain,
\item or a factorial Dedekind domain.
\end{enumerate}
Then $g\in Q(S)$ is of geometric origin if and only if  $\ob(g)\in H^2(S,\GG_m)$ vanishes.
\end{corollary}

\proof
According to  Hilbert 90, every $\GL_n$-torsor over a scheme is trivial in the Zariski topology.
If $R$ is semilocal,   the set  $H^1(S,\GL_n)$ is thus a singleton, and the assumption of the Theorem 
obviously holds. 

If $R$ is a Dedekind ring, the structure theory of locally free modules of finite rank
gives that the determinant map $H^1(S,\GL_n)\ra H^1(S,\GG_m)$ is bijective. If $\Pic(S)=0$, in other words $R$ is factorial,
the set $H^1(S,\GL_n)$ must be a singleton.
By the Quillen--Suslin Theorem (\cite{Quillen 1976} and  \cite{Suslin 1976}),  the same is true if  $R=k[T_1,\ldots,T_r]$ or $R=D[T_1,\ldots,T_r]$,
where $k$ is a field and $D$ is a principal ideal domain.
\qed

%===========================================================
\section{Spaces of Schur representations}
\mylabel{Spaces representations}

Throughout, an \emph{associative ring} is an abelian group $\Lambda$, endowed with 
a bilinear multiplication that possesses a unit element and satisfies the axiom of associativity.
If additionally    the axiom
of commutativity   holds it is  simply called \emph{ring}. An \emph{associative  algebra} over a ring $R$ is an associative ring $\Lambda$,
together with a homomorphism $\varphi:R\ra \Lambda$ whose image $\varphi(R)$ belongs to the center
$Z(\Lambda)$. Then for each ring homomorphism $R\ra A$, the tensor product $\Lambda\otimes_RA$ carries in a canonical way the
structure of an $A$-algebra.

Let $R$ be a ground ring,
and $\Lambda$ and $\Lambda'$ be two associative $R$-algebras.
Important special cases arise when $\Lambda=R[\Gamma]$ is a group algebra of some group $\Gamma$,
or the path algebra $\Lambda=\pathalg{R}$ of some quiver $\quiver$, and $\Lambda'$ is a matrix algebra or more generally an   Azumaya algebra.
For each $R$-algebra $A$ we   can form the set
$$
F(A)=F_{\Lambda/R}^{\Lambda'}(A)=\Hom_{\text{$A$-Alg}}(\Lambda\otimes_R A,\Lambda'\otimes_R A).
$$
The construction is functorial in $A$, and defines a set-valued contravariant functor
on the category $(\Aff/R)$ of affine $R$-schemes, which is opposite to the category of   $R$-algebras.
The elements $\rho\in F(A)$ are called  \emph{$A$-valued representations of $\Lambda$ in $\Lambda'$}.

\begin{proposition}
\mylabel{sheaf axiom}
The above contravariant functor $F:(\Aff/R)\ra (\Set)$ satisfies the sheaf axiom with respect to the fpqc topology.
\end{proposition}

\proof
Let $A$ be an $R$-algebra, $A\subset A_0$ be a faithfully flat ring extension, and   $A_1=A_0\otimes_AA_0$.
We have to check that 
$$
\begin{tikzcd} 
F(A) \ar[r]	& F(A_0)\ar[r, shift left=.5ex] \ar[r, shift right = .5ex]	& F(A_1) 
\end{tikzcd}
$$
is an equalizer diagram of sets. In other words, the arrow on the left is injective, and its image comprises
the $\rho_0\in F(A_0)$ whose two images    in $F(A_1)$ coincide.
Consider the larger functor $G(A)=\Hom_{\text{$A$-Mod}}(\Lambda\otimes_R A,\Lambda'\otimes_R A)$ given by linear maps,
not necessarily preserving multiplication and unit element. By fpqc descent, the above diagram with $G$ instead of $F$ is
an equalizer (\cite{SGA 1}, Expos\'e VIII, Corollary 1.2). It follows that $F(A)\ra F(A_0)$ is injective, and each $\rho_0\in F(A_0)$ whose images
in $F(A_1)$ coincide  descents to some $\rho\in G(A)$.
It remains to check that $\rho(xy)=\rho(x)\rho(y)$ for all $x,y\in \Lambda\otimes_RA$, and $\rho(1)=1$.
This follows from the commutative diagram
$$
\begin{CD}
\Lambda\otimes_R A		@>\iota>>	\Lambda\otimes_RA_0\\
@V\rho VV		@VV\rho_0V\\
\Lambda'\otimes_R A		@>>\iota'>	\Lambda'\otimes_RA_0,
\end{CD}
$$
because the canonical map  $\iota'$ is injective, and $\iota$, $\rho_0$ and $\iota'$ respect multiplications and unit element.
\qed

\medskip
Suppose from now on that the structure map $R\ra \Lambda'$ is locally a direct summand of $R$-modules; this ensures that
all induced maps $A\ra \Lambda'\otimes A$ remain injective. We then consider the subsets
$$
F^0(A)\subset F(A)
$$
comprising  the   homomorphisms $\rho:\Lambda\otimes_RA\ra\Lambda'\otimes_RA$
such that for each prime ideal $\primid\subset A$, the only elements in $\Lambda'\otimes_R\kappa(\primid)$ that
commute with all $\rho(x)$, $x\in\Lambda\otimes_R\kappa(\primid)$ are the members of $R\otimes_R\kappa(\primid)=\kappa(\primid)$,
in other words, the scalars.
Clearly, the formation of $F^0(A)\subset F(A)$ is functorial in $A$, and thus defines a subfunctor $F^0\subset F$.
The sheaf axiom with respect to the fpqc topology   holds for $F^0$, because it holds for $F$,
and the formation of commutants commutes with field extensions.
The elements $\rho\in F^0(A)$ are called $A$-valued \emph{Schur representations} of $\Lambda$ in $\Lambda'$.

\begin{proposition}
\mylabel{relatively representable}
Suppose   $\Lambda'$ is locally free  of finite rank  as $R$-module.
Then the inclusion $F^0\subset F$ is relatively representable by open embeddings.
\end{proposition}

\proof
Let $V=\Spec(A)$ be an affine $R$-scheme,   and $\rho:\Lambda\otimes_RA\ra\Lambda'\otimes_RA$ be some homomorphism.
By the Yoneda Lemma, we may  regard it  as a natural transformation
$\rho:V\ra F$, and have to check that the ensuing fiber product $F^0\times_FV$ is representable by a scheme,
and that the projection to $V$ is an open  embedding. Clearly, the projection $F^0\times_FV\ra V$ is a monomorphism.

Without loss of generality we may assume $A=R$, and that the $R$-module  $\Lambda'$ is free of rank $r\geq 1$. 
Let $U\subset V$ be the set of points   $a\in V$ where the canonical morphism $\Spec \kappa(a)\ra V$
factors over $F^0\times_FV$. Suppose for the moment that this is an open set, and thus
defines an open subscheme. By the very definition of $F^0$ and $U$, the two monomorphisms
$F^0\times_FV\ra V$ and $U\ra V$ factor over each other, hence the former is represented by the latter.

It remains to verify that the subset $U$ is indeed open. Fix a point $a\in U$, in other words,
$\kappa(a)\subset\Lambda'\otimes\kappa(a)$ is the commutant for $\rho\otimes\kappa(a)$. Using that  
the vector space  $\Lambda'\otimes\kappa(a)$ is finite-dimensional, we find  finitely many $g_1,\ldots,g_n\in \Lambda$ such that $\kappa(a)$ is the commutant
for the $\rho(g_i)\otimes 1$, $1\leq i\leq n$.
Consider the $R$-linear map
$$
\Psi:\Lambda'\lra \bigoplus_{i=1}^n\Lambda',\quad f\longmapsto (f\rho(g_i)-\rho(g_i)f)_{1\leq i\leq n}
$$
between free $R$-modules of finite rank. For each point $v\in V$, the kernel for $\Psi\otimes\kappa(v)$ 
contains $\kappa(v)$, hence $\rank(\Psi\otimes\kappa(v))\leq r-1$.
For $v=a$ this becomes an equality.
Viewing $\Psi$ as an $rn\times r$-matrix, we see that     some $(r-1)$-minor $h\in R$ 
does not vanish in $\kappa(a)$.
Replacing $R$ by the localization $R[h^{-1}]$ we may assume that the minor is a unit. This ensures that
the function $v\mapsto \rank(\Psi\otimes\kappa(v))$ takes constant value $r-1$.  
By \cite{EGA I}, Chapter 0, Proposition 5.5.4 the image of $\Psi$ is locally free of rank $r-1$, hence the kernel is invertible.
For each point $v\in V$ the unit element $1\in\Lambda'$ generates $\Kernel(\Psi)\otimes\kappa(v)$,
so the inclusion $R\subset \Kernel(\Psi)$ is an equality. It follows that $\kappa(v)\subset\Lambda\otimes\kappa(v)$
is the commutant for the $\rho(g_i)\otimes 1$, and thus for $\rho\otimes\kappa(v)$, for all $v\in V$.
This shows $U=V$.
\qed

\medskip
We also have   an absolute representability statement:

\begin{proposition}
\mylabel{schur representable}
Suppose  $\Lambda'$ is locally free of finite rank as  $R$-module.
Then   $F:(\Aff/R)\ra (\Set)$ is representable by an affine scheme.
It is of finite type provided that $\Lambda$ is finitely generated.
\end{proposition}

\proof
Suppose first that the associative algebra $\Lambda=R\langle T_i\rangle_{i\in I}$ is free.
One  easily checks that $F$ is represented by the spectrum of
$$
B=\bigotimes_{i\in I}\Sym^\bullet(\Lambda'^\vee)=\Sym^\bullet(\bigoplus_{i\in I}\Lambda'^\vee),
$$
by using various universal properties and   biduality $\Lambda'=\Lambda'^{\vee\vee}$.

For the general case, express  the associative algebra $\Lambda$ in terms of generators  $g_i\in\Lambda$, $i\in I$
and relations $r_j\in R\langle T_i\rangle_{i\in I}$, $j\in J$.
The preceding paragraph gives a monomorphism $F\subset \Spec(B)$, and we have to check that this is relatively
representable by closed embeddings. This is a local problem, so we may assume that $\Lambda'$ admits
a basis $e_1,\ldots,e_r$ as $R$-module. For each $P\in R\langle T_i\rangle_{i\in I}$, the equation
$$
P(\ldots, x_i ,\ldots ) = \sum_{k=1}^r P_k(\ldots, x_i ,\ldots ) \cdot e_k
$$
inside $\Lambda'$ defines  polynomial maps $P_k:\bigoplus_{i\in I}\Lambda'\ra R$.  These maps stem from 
elements $P_k\in \Sym^\bullet(\bigoplus_{i\in I}\Lambda'^\vee)=B$, because multiplication in $\Lambda'$ is bilinear.
One easily checks that $F\subset\Spec(B)$ is the closed subscheme
defined by the ideal generated by the $P_1,\ldots,P_r$, where $P=r_j$, $j\in J$ ranges over the relations. 
\qed

\medskip
Suppose now that $\Lambda'$ is an \emph{Azumaya algebra}  of   degree $n\geq 1$, that is, a twisted form of $\Mat_n(R)$.
Note that $\Lambda'=\Mat_n(R)$ is indeed an important special case.
Write $U_{\Lambda'/R}$ and $\Aut_{\Lambda'/R}$ for the group-valued sheaves on $(\Aff/R)$ defined by
$$
U_{\Lambda'/R}(A) = (\Lambda'\otimes_RA)^\times \quadand \Aut_{\Lambda'/R}(A) = \Aut(\Lambda'\otimes_RA).
$$
These are twisted forms of $\GL_n$ and $\PGL_n$, confer the discussion in \cite{Schroeer; Tziolas 2023}, Lemma 3.1. 
In particular, $U_{\Lambda'/R}$ and $\Aut_{\Lambda'/R}$ are smooth affine group schemes
of finite type. The     conjugacy map defines an exact sequence
$$
1\lra \GG_m\lra U_{\Lambda'/R}\stackrel{\conj}{\lra} \Aut_{\Lambda'/R}\lra 1,
$$
where the inclusion on the left comes from the structure inclusion $R\subset\Lambda'$.

The group scheme $G=\Aut_{\Lambda'/R}$ acts on the sheaf $F:(\Aff/R)\ra (\Set)$ via composition $(g,\rho)\mapsto g\circ\rho$.
Now let $\tilde{\Lambda}'$ be a twisted form of $\Lambda'$, that is, another Azumaya algebra of the same degree $n$. 
Write $P$ for the corresponding $G$-torsor,
and let $\tilde{F}:(\Aff/R)\ra (\Set)$ be the sheaf defined with $\Lambda$ and $\tilde{\Lambda}'$.

\begin{proposition}
\mylabel{schur subfunctor}
Notation as above. Then $F^0\subset F$ is $G$-stable, the induced $G$-action on $F^0$ is free, and we have  canonical identifications
$$
P\wedge^G F = \tilde{F}\quadand P\wedge^G F^0 = \tilde{F}^0
$$
of contravariant functors on $(\Aff/R)$.  
\end{proposition}

\proof
The open subfunctor $F^0$ is clearly $G$-stable. Suppose we have an $R$-algebra $A$, and elements $g\in G(A)$ and $\rho\in F^0(A)$
with $g\circ \rho=\rho$.
Choose some fpqc extension $A\subset A_0$ so that $\Lambda'\otimes_RA_0\simeq\Mat_n(R_0)$.
Then $g_0=g\otimes A_0$ becomes an element in  $\PGL_n(A_0)$. Replacing $A_0$ by some further fpqc extension,  we may assume
that it stems from some $S\in \GL_n(A_0)$. Setting $\rho_0=\rho\otimes A_0$, we obtain  $g_0\circ \rho_0=S\cdot\rho_0\cdot S^{-1}$.
In other words, $S\cdot\rho_0=\rho_0\cdot S$. Using that $\rho$ is a Schur representation, we see that $S$ must be a scalar matrix.
Thus $g_0$ and hence also $g$ are trivial. Consequently,  the $G$-action on $F^0$ is free.

For the last assertion we first assume  that there is an isomorphism $\psi:\Lambda'\ra\tilde{\Lambda}'$ of Azumaya algebras.
This corresponds to an equivariant map $\psi: G\ra P$, and the arrows in
$$
P\wedge^GF\stackrel{\psi}{\longleftarrow}G\wedge^GF = F\stackrel{\psi}{\lra} \tilde{F}
$$
yield the desired canonical identification $P\wedge^G F = \tilde{F}$.

In the general case, we find an fpqc extension $R\subset R_0$ and   an isomorphism $\psi_0:\Lambda'\otimes_RR_0\ra\tilde{\Lambda}'\otimes_RR_0$,
with corresponding  equivariant $\psi_0:G\otimes_RR_0\ra P\otimes_RR_0$.
This yields $P\wedge^G F = \tilde{F}$ if we restrict the contravariant sheaves along the forgetful functor $(\Aff/R_0)\ra (\Aff/R)$.
Since both $P\wedge^G F$ and $\tilde{F}$ satisfy the sheaf axiom with respect to the fpqc topology, and since every fpqc extension $A\subset A'$ 
of $R$-algebras can be refined to an fpqc covering $A\otimes_RR_0\subset A'\otimes_RR_0$ of $R_0$-algebras, the identification $P\wedge^G F = \tilde{F}$ 
already holds over $(\Aff/R)$. Since the $G$-action on $F$ stabilizes $F^0$, we get an induced identification $P\wedge^G F^0 = \tilde{F}^0$.
\qed

\medskip
Let us summarize our findings:

\begin{theorem}
\mylabel{space  schur representations}
Suppose $\Lambda$ is finitely presented as an  associative $R$-algebra,
and $\Lambda^\azu$ is an Azumaya algebra of degree $n\geq 1$. Then the functor $X=X^{\Lambda^\azu}_{\Lambda/R}$ of Schur representations
is   representable by a  quasiaffine scheme  of finite presentation. Moreover, the group scheme $H=\Aut_{\Lambda^\azu/R}$ acts freely,
and the quotient $Q=X/H^\op$ is an algebraic space that is  of finite presentation.
\end{theorem}

\proof
By Proposition \ref{relatively representable} and \ref{schur representable}, 
the functor $X$ is representable by some open set in some affine scheme of finite type. 
Obviously, the structure morphism $H\ra S$ is flat and locally of finite presentation. Since the $H$-action on $X$ is free,
the quotient $Q=X/H^\op$ exists as an algebraic space (see for example \cite{Laurent; Schroeer 2024}, Lemma 1.1). 
It follows from \cite{EGA IVb}, Proposition 2.7.1 that $Q$ is locally of finite type.
Moreover, $Q$ is of finite presentation if  this holds for $X$.

It remains to check that $X$ is quasicompact  and  locally of finite presentation. 
In light of loc.\ cit., we may replace $R$ by some fppf extension,
and assume that $\Lambda^\azu=\Mat_n(R)$. Since the associative algebra $\Lambda$ involves only finitely many structure constants,
we may furthermore assume that the ground ring $R$ is finitely generated over $\ZZ$. Using that $X$ is open in
some affine scheme of finite type, we see that it is noetherian. 
\qed

\medskip
To close this section, let us relate  Schur representations to other notions from representation theory,
for simplicity when $R=K$ is a field. Let $M$ be a $\Lambda$-module that is finite-dimensional as $K$-vector space,
and $\Endomor=\End_\Lambda(M)$ be its endomorphism ring.  One says that $M$ is \emph{geometrically simple}  if for all field extensions $K\subset L$,
the base change $M_L=M\otimes_KL $ is simple as representation of $\Lambda_L=\Lambda\otimes_K L$. 
The notion of \emph{geometrically indecomposable}
is formed in a similar way. The alternative terms \emph{absolutely simple} and \emph{absolutely indecomposable} are also frequently in use.
Recall that the associative ring $\Endomor$ is  
\emph{local} if the non-units form a left ideal $\maxid$. This is indeed maximal and two-sided, giving the residue skew field $\kappa=\Endomor/\maxid$.

\begin{lemma}
\mylabel{simple and indecomposable}
In the above situation, the following holds:
\begin{enumerate}
\item If the module $M$ is simple, then  $\Endomor$ is a skew field. 
\item The module $M$ is indecomposable if and only if    $\Endomor$ is  local.
\item If $M$ is geometrically simple, then $M$ is Schur.
\end{enumerate}
Moreover, if  $M$ is indecomposable and $L$ is   purely inseparable,  then the $\Lambda_L$-module $M_L$ is indecomposable.
\end{lemma}
 
\proof
The first assertion  is Schur's Lemma. Statement (ii) can be found in \cite{Curtis; Reiner 1990}, Proposition 6.10.
For (iii) choose a separably closed extension $L$. Since  $M_L$ is simple and the Brauer group $\Br(L)$ vanishes, the   inclusion $L\subset\Endomor_L$ must
be an equality. Hence the same holds for $K\subset\Endomor$, and $M$ is Schur.

Finally, let $M$ be indecomposable and $L$   purely inseparable. To check that $M_L$ is indecomposable it suffices
to treat the case that $L\subset K^{1/p}$, where $p>0$ is the characteristic. Let  $f,g\in \Endomor_L$ be two non-units.
Then $f^p,g^p$ are non-units that belong to $\Endomor\subset\Endomor_L$. They vanish in the residue skew field $\kappa$, so the
same holds for $f,g$. It follows that $f+g$ vanishes in the residue class ring $\kappa_L$ of $\Endomor_L$,
so $f+g\in\Endomor_L$ is a non-unit. Consequently, the non-units   form a left ideal.
\qed

\medskip
Note that Kraft and Riedtmann (\cite{Kraft; Riedtmann 1986}, Theorem in 2.6) showed that a quiver representation $M$ is 
Schur if and only if the point $M$ in the representation space  admits an open neighborhood where all points are geometrically
indecomposable.

Also note that moduli spaces of representations are traditionally formed in the realm of \emph{Geometric Invariant Theory} \cite{Mumford; Fogarty; Kirwan 1993}, which
relies on the choice of   \emph{stability conditions}. The latter ensure, by design, that the \emph{geometrically stable objects} are Schur representations. 
In fact the geometrically stable objects for any given stability condition form an open subset of $X$. In turn, all  possible  GIT quotients are 
simultaneously contained   as schematic open subsets in our algebraic space $Q=X/H^\op$. However, the semistable objects, which frequently appear on the boundary in GIT quotients, 
often fail to be Schur and are thus beyond the scope of this paper.

%===========================================================
\section{Quiver representations}
\mylabel{Quiver representations}

\newcommand{\source}{\operatorname{source}}
\newcommand{\target}{\operatorname{target}}
Our initial motivation for this paper was to understand twisting of   quiver representations,
as studied by Hoskins and Schaffhauser \cite{Hoskins; Schaffhauser 2020}.
In this   section we indeed take a closer look at   representations of quivers in Azumaya algebras,
as a special case of the  general theory developed so far, and show that the resulting space of Schur representations
embeds into a product of Grassmann varieties and vector bundles.

Recall that a \emph{quiver}  is a finite  directed graph $\quiver$, with loops and  multiple arrows allowed.
Formally, it comprises a set of vertices $\quiver_0$ and a set of arrows $\quiver_1$, together with   source and target maps $\quiver_1\rightrightarrows \quiver_0$.
An arrow $\alpha$ starting at a  vertex $i$ and ending at a vertex $j$   thus has  $i=\source(\alpha)$ and $j=\target(\alpha)$. 
By abuse of notation, we write   $i,\alpha\in \quiver$ instead of $i\in\quiver_0 , \alpha\in\quiver_1$ 
to denote both vertices and arrows from the quiver. No confusion is possible, since we always use Latin letters for vertices,
and Greek letters for arrows. 
A  \emph{quiver representation}  over a field $k$ is a collection of   vector spaces $V_i$ for the vertices $i\in \quiver$,
together with linear maps $f_\alpha: V_{\source(\alpha)}\ra V_{\target(\alpha)}$ for the arrows $\alpha\in \quiver$.
We usually write this datum in the form  $(V_i\mid f_\alpha)_{i,\alpha\in \quiver}$. With the obvious notion of homomorphisms,
the quiver representations form a $k$-linear abelian category.
 
Now fix  a ground ring $R$  and write  $\Lambda=\pathalg{R}$ for the resulting \emph{path algebra},
which is the associative $R$-algebra generated by formal symbols $e_i$ and $f_\alpha$, 
for the vertices and  arrows $i,\alpha\in \quiver$, subject to the relations
\begin{equation}
\label{path algebra relations}
e_ie_j=\delta_{i,j}e_i\quadand f_\alpha e_i=\delta_{\source(\alpha),i}f_\alpha\quadand  e_jf_\beta=\delta_{j,\target(\beta)}f_\beta,
\end{equation}
where $\delta_{i,j}$ denotes the Kronecker delta.  
From the relations  one immediately sees
that multiplication with  $\sum_{i\in \quiver}e_i$, from either the left or   the right,  fixes each generator, and thus
 $\sum_{i\in \quiver}e_i=1$. One also  says that the $e_i\in\Lambda$, $i\in\quiver$ form a \emph{partition of unity into orthogonal idempotents}. 
Recall that quiver representations over $R$ can be identified with left modules over the path algebra $\pathalg{R}$ 
(compare \cite{Kirillov 2016}, Theorem 1.7). 

We will now generalize this correspondence to representations of $\pathalg{R}$ in an Azumaya algebra $\Lambda^\azu$.
Let $U$ be a scheme and $\shA$ be an associative $\O_U$-algebra whose underlying $\O_U$-module is quasicoherent. 
We call a datum of the form $(\shV_i\mid f_\alpha)_{i,\alpha\in \quiver}$, 
where $\shV_i$ is a sheaf of right $\shA$-modules and 
$f_\alpha: \shV_{\source(\alpha)}\ra \shV_{\target(\alpha)}$ is a morphism of $\shA$-modules, 
a \emph{quiver representation of $\quiver$ in right $\shA$-modules}. 
We will discuss at the end of this Section that this indeed generalizes the ordinary notion of quiver representations.

Now let $\rho: \O_U[\quiver]=\Lambda\otimes_R\O_U\ra \shA$
be a homomorphism  of associative $\O_U$-algebras.
The   subsheaves $$e_i \cdot \shA=\rho(e_i\otimes 1)\cdot \shA\subset\shA$$
yield a direct sum decomposition 
$\bigoplus_{i\in \quiver}( e_i\cdot\shA)=\shA$ of right $\shA$-modules, in light of $\sum_{i\in \quiver}e_i=1$. 
By the relations \eqref{path algebra relations}, each arrow $\alpha\in \quiver$ defines 
$\shA$-linear maps between the summands via left multiplication with $\rho(f_\alpha\otimes 1)$. 
By   abuse of notation these maps are  denoted by 
$f_\alpha: e_i\cdot \shA\ra  e_j\cdot \shA$.
All of them vanish, except for $i=\source(\alpha)$ and $j=\target(\alpha)$.
Hence, $(e_i\cdot\shA\mid f_\alpha)_{i,\alpha\in \quiver}$ is a quiver representation of $\quiver$ in right $\shA$-modules and we
will see in the proof of Proposition \ref{embedding for quiver} below that
$\rho:\O_U[\quiver]\ra \shA$ is already determined by  $(e_i\cdot\shA\mid f_\alpha)_{i,\alpha\in \quiver}$.
Note that if $\shA$  is  locally free of finite rank as an $\O_U$-module, the same holds
for the summands $e_i\cdot\shA$.

Fix an Azumaya algebra $\Lambda^\azu$ of degree $n\geq 1$, and set $H=\Aut_{\Lambda^\azu/R}$.
We now examine the quasiaffine scheme $X=X^{\Lambda^\azu}_{\pathalg{R}/R}$
of Schur representations  in the Azumaya algebra  $\Lambda^\azu$.  
Our observation here is  that $X$ can also be described 
in terms of Grassmann varieties and vector bundles:
 Given a projective $R$-module $E$ of finite rank $r\geq 0$, we write  
$$
\AA^1\otimes_RE =\Spec\Sym^\bullet(E^\vee)\quadand \Grass_{E^\vee/R}^d 
$$
for the corresponding vector bundle of rank $r$ and Grassmann varieties of relative dimension $d(r-d)$. 
The elements in the bidual $E=E^{\vee\vee}$ and the locally free quotients   $\psi:E^\vee\ra M$ of rank $d$ 
give the respective sets of $R$-valued points.
The latter can also be seen as the locally direct summands $N\subset E$ of corank $d$, via  $N=M^\vee$. 
By the discussion above, we have a canonical morphism
\begin{equation}
\label{grassmann varieties}
X=X^{\Lambda^\azu}_{\pathalg{R}/R}
\subset F^{\Lambda^\azu}_{\pathalg{R}/R} \lra 
\dot{\bigcup_m} \left(\prod_{i\in \quiver} \Grass^{n^2-m_i}_{(\Lambda^\azu)^\vee/R} \times \prod_{\alpha\in\quiver}
\AA^1\otimes_R\Lambda^\azu \right),
\end{equation} 
where the disjoint union is indexed by tuples $m=(m_i)_{i\in \quiver}$ of natural numbers  
subject to $\sum_{i\in \quiver}m_i=n^2$. 
It sends  an $A$-valued path algebra representation $\rho:\pathalg{A}\ra\Lambda^\azu\otimes_RA$  
to the  quiver representation 
%datum 
$(e_i\cdot \Lambda^\azu\otimes_RA\mid f_\alpha)_{i,\alpha\in \quiver}$. 
%Note that the 
 The group scheme
$H=\Aut_{\Lambda^\azu/R}$ acts in a canonical way on $X$, and also on the disjoint union on the right.
Note that by our choice of   conventions, the action is in both cases indeed from the left.
 
\begin{proposition}
\mylabel{embedding for quiver}
The above morphism is  $H$-equivariant, and an   embedding of quasiprojective schemes.
\end{proposition}

\proof
Write $Y$ for the right-hand side of \eqref{grassmann varieties}. The group elements $\sigma\in H(A)$ act via
$$
\sigma\cdot\varphi = \sigma\circ\varphi\quadand
\sigma\cdot (e_i\cdot \Lambda^\azu\otimes_RA\mid f_\alpha)=
 (\sigma(e_i \cdot \Lambda^\azu\otimes_RA)\mid \sigma\circ f_\alpha\circ \sigma^{-1}),
$$
so equivariance of $X\ra Y$ is clear. We now actually prove that this morphism is relatively representable by embeddings.
Note that this gives, for   quiver representations,   
a proof for representability   independent of Theorem \ref{space schur representations}.

Let us start by checking that $X\ra Y$ is a monomorphism.
Suppose we have 
 an $A$-valued representation $\rho:\pathalg{A}\ra \Lambda^\azu\otimes_RA$ of the path algebra,
and set $\ideala_i=e_i\cdot\Lambda^\azu\otimes_RA$. 
We need to show that $\rho$ can be reconstructed  
from the datum $(\ideala_i\mid f_\alpha)_{i,\alpha\in \quiver}$, and for this it suffices to treat the case $A=R$.
As explained in \cite{A 8}, Chapter VIII, \S 8, No.\ 4 
the decomposition  $\Lambda^\azu=\bigoplus_{i\in\quiver}\ideala_i$ into right ideals corresponds to a
partition of unity into orthogonal idempotents $e_i\in\Lambda^\azu$, $i\in \quiver$. 
Together with the $f_\alpha\in \Lambda^\azu$, $\alpha\in\quiver$ we recover the path algebra representation
$\rho:\pathalg{R}\ra\Lambda^\azu$.
% Since the natural map $e_j\cdot\Lambda^\azu\cdot e_i \to \Hom_{\Lambda^\azu}(e_i\cdot\Lambda^\azu,e_j\cdot\Lambda^\azu)$ 
% given by left multiplication is bijective,
% we  recover the representation $\rho:\pathalg{R}\ra\Lambda^\azu$. 

Given a morphism $S\ra Y$ from some affine scheme $S=\Spec(A)$,  it remains to verify
that $X'=X\times_YS$ is representable by a scheme, and that the projection to $S$  is an embedding. In other words, it
factors as  a closed embedding inside  some open set.
Again it suffices to treat the case $A=R$.  Let $(\ideala_i\mid f_\alpha)_{i,\alpha\in\quiver}$ be the datum
defining  $S\ra Y$.  
So $\ideala_i\subseteq \Lambda^\azu$ are locally free $R$-submodules and $f_\alpha\in\Lambda^\azu$ are arbitrary elements.

The cokernel $M$ for the canonical mapping $\bigoplus_{i\in\quiver}\ideala_i\ra\Lambda^\azu$ 
is of finite presentation, so its support defines a closed set $Z\subset S$, and $X'\ra S$ factors
over the complementary open set. Replacing $R$ by  suitable localizations, we thus may
assume that $\bigoplus_{i\in\quiver}\ideala_i\ra\Lambda^\azu$ is surjective. Since both $R$-modules are locally free
of the same rank, the map is actually bijective.

Next we fix a vertex $i\in\quiver$, and choose    algebra generators $g_1,\ldots,g_r\in\Lambda^\azu$.
The $R$-submodule $\ideala_i\subset\Lambda^\azu$
is a  right  ideal if and only if the cokernels $M_{ij}$ for the right multiplications
$\ideala_i\stackrel{g_j}{\ra} \Lambda^\azu/\ideala_i$ are  zero.
The support of $M_{ij}$ defines a closed set and $X'\ra S$ factors over the complementary open set.
As in the preceding paragraph, we may replace $R$ by suitable localizations, and assume that the $\ideala_i$ are right ideals.

Now  we fix an arrow $\alpha\in\quiver$, and write $i,j\in\quiver$ for its source and target. 
The right multiplication $f_\alpha:\Lambda^\azu\ra\Lambda^\azu$ sends $\ideala_i$ to $\ideala_j$ and vanishes on $\ideala_r$, $r\neq i$
if and only if the
induced maps   $\ideala_i\ra \Lambda^\azu/\ideala_j$ and $\ideala_r\ra\Lambda^\azu$ vanish.
These are  closed conditions.  

This shows  that $F=F^{\Lambda^\azu}_{\pathalg{R}/R}\ra Y$ is relatively representable by embeddings. 
Since $X=F^0\subset F$ is relatively representable by open embeddings by Proposition \ref{relatively representable}, 
we obtain that $X\ra Y$ is relatively representable by embeddings as well.
\qed

\medskip
To close this section, let us discuss the special case of matrix algebras, and verify
that the quiver representations in right modules   over $\Lambda^\azu=\Mat_n(R)$ corresponds to the classical notation
of quiver representations    over $R$. Recall that the functors
$$
M\longmapsto M\otimes_{\Mat_n(R)} \Mat_{n\times 1}(R)\quadand
N\longmapsto N\otimes_R\Mat_{1\times n}(R)
$$
are   quasi-inverse equivalences between the category of right modules $M$ over the matrix algebra $\Mat_n(R)$ and modules $N$ over
the commutative ring $R$
(Morita equivalence, compare also the result of Watts \cite{Watts 1960}).
Set $V=\Mat_{n\times 1}(R)$. Under the above equivalences, the 
decomposition $\Mat_n(R)=\bigoplus_{i\in\Delta} \ideala_i$ into right ideals corresponds the   decomposition
$V=\bigoplus_{i\in\Delta} V_i$ into submodules $V_i$, where the summands  are locally free of some rank $v_i$, subject
to $\sum v_i=n$. From $\ideala_i=V_i\otimes_R\Mat_{1\times n}(R)$ we see $v_i=m_i/n$.

In the theory of quiver representations, the tuple $(v_i)_{i\in \Delta}$ is called \emph{dimension vector}.
Analogously to Proposition \ref{embedding for quiver} one obtains an embedding
\begin{equation}
\label{equ_ordinary_embedding}
X  =X^{\Mat_n(R)}_{\pathalg{R}/R}\lra 
\dot{\bigcup_v} \left(\prod_{i\in \quiver} \Grass^{n-v_i}_{\Mat_{n\times 1}(R)^\vee/R} \times 
\prod_{\alpha\in\quiver}\AA^1\otimes_R\Mat_n(R) \right),
\end{equation}
where the disjoint union now runs over all tuples $v=(v_i)_{i\in\Delta}$ subject to $\sum v_i=n$.
In the theory of quiver representations, it is customary to fix the dimension vector $v$,  and
a decomposition $\Mat_{n\times 1}(R)=\bigoplus_{i\in \Delta} V_i$ with $\rank(V_i)=v_i$.
The corresponding closed subscheme  $X'\subset X$ becomes an open subscheme of the vector  bundle
$\prod_{\alpha\in\quiver} \AA^1\otimes_R\Hom (V_{\source(\alpha)}, V_{\target(\alpha)})$, which can be seen as a subbundle in
$\prod_{\alpha\in\quiver}\AA^1\otimes_R\Mat_n(R)$. Note that this also shows  that $X$ is smooth,
a fact that could also be proven by using that path algebras of quivers are formal smooth,
in the sense of associative algebras.

%===========================================================
\section{The gerbe of splittings for an Azumaya algebra}
\mylabel{Gerbe splittings}

In this section we develop some   facts on splittings of Azumaya algebras in the context of stacks. The material is well-known
(confer the work of Lieblich \cite{Lieblich 2004}, \cite{Lieblich 2007}, \cite{Lieblich 2013}),
but for our purposes 
we need it  in a form that makes the involved fibered categories, the resulting topoi of sheaves, and the continuous maps between them explicit.

Fix a ground ring $R$, and some Azumaya algebra $\Lambda^\azu$ of degree $n\geq 1$.
For any scheme $V$, let us  write   $V'$ for the  \emph{$\GG_{m,V}$-gerbe of splittings} for the Azumaya algebra $\Lambda^\azu\otimes_R\O_V$.
This is the category
$$
V'=\{(U,h, \shF,\psi)\mid \text{ $\psi:\uEnd(\shF)\ra \Lambda^\azu\otimes_R\O_U$}\}
$$
whose objects are  quadruples $(U,h, \shF,\psi)$ where   $U$ is an affine scheme,
$h:U\ra V$ is a morphism of   schemes, $\shF$ is a locally free sheaf of rank $n$ over $U$,
and $\psi:\uEnd(\shF)\ra \Lambda^\azu\otimes_R\O_U$ is an isomorphism of sheaves of associative algebras.  

Arrows $(U_1,h_1, \shF_1,\psi_1)\ra (U_2,h_2, \shF_2,\psi_2)$ in the category $V'$ are pairs $(g,\Psi)$ where 
$g:U_1\ra U_2$ is a morphism of schemes with $h_2\circ g=h_1$, and $\Psi:\shF_2\ra g_*(\shF_1)$ is a homomorphism of quasicoherent sheaves
such that the adjoint $\Phi:g^*(\shF_2)\ra\shF_1$ is an isomorphism making the diagram
$$
\begin{tikzcd} 
g^*(\shF_2)^\vee\otimes g^*(\shF_2)\ar[rr,"\Phi^{\vee-1}\otimes\Phi"]\ar[dr,"g^*(\psi_2)"']	& & \shF_1^\vee\otimes\shF_1\ar[dl,"\psi_1"]\\
					& \Lambda^\azu\otimes_R\O_{U_1}
\end{tikzcd}
$$
commutative.  Here we employ the identifications $\shF_i^\vee\otimes\shF_i=\End(\shF_i)$ discussed in the proof for Lemma \ref{equivariant} below.
Using tensor products and fppf descent for quasicoherent sheaves, one easily checks that the forgetful functor
$$
V'\lra (\Aff/V),\quad (U,h, \shF,\psi)\longmapsto (U,h)
$$
endows $V'$ with the structure of a fibered category satisfying the stack axioms. The Skolem--Noether Theorem ensures
that the fiber categories are groupoids. Note that despite our notation, $V'$ is a stack rather than a scheme. 

The category $V'$ carries a Grothendieck topology, where 
$(U_\lambda,h_\lambda, \shF_\lambda,\psi_\lambda)_{\lambda\in L}\ra (U,h, \shF,\psi)$ is a covering family if  $(U_\lambda\ra U)_{\lambda\in L}$
is an fppf covering of   affine schemes. In this way we regard $V'$ as a site, with the ensuing notion of sheaves $F'$ on $V'$.
We remark in passing that there is    the Zariski topology, the \'etale topology, and the fpqc topology as well.
The site $V'$ comes  with a \emph{structure sheaf} $\O_{V'}$ and a  \emph{tautological sheaf} $\shF^\taut_{V'}$, defined 
for objects  $\tilde{U}=(U,h, \shF,\psi)$ by 
\begin{equation}
\label{sheaves on gerbe}
\Gamma(\tilde{U},\O_{V'}) =\Gamma(U,\O_U)\quadand \Gamma(\tilde{U},\shF^\taut_{V'}) =\Gamma(U,\shF),
\end{equation}
with obvious restriction maps.
We also have a  \emph{tautological splitting}, which is the isomorphism
$$
\psi^\taut_{V'}:\uEnd(\shF^\taut_{V'})\lra\Lambda^\azu\otimes_R\O_{V'}
$$
constructed as follows:  For each endomorphism $\alpha$ of $\Gamma(\tilde{U},\shF^\taut_{V'})=\Gamma(U,\shF)$ that is linear with respect to 
$\Gamma(\tilde{U},\O_{V'})=\Gamma(U,\O_U)$,
we set $\psi^\taut_{V'}(\alpha)=\psi(\alpha)$ as elements in  
$$
\Gamma(U,\Lambda^\azu\otimes_R\O_U) = \Lambda^\azu\otimes_R\Gamma(U,\O_U)= \Gamma(\tilde{U},\Lambda^\azu\otimes_R\O_{V'}).
$$
Obviously, this is compatible with restrictions.

Given a morphism $f:V_1\ra V_2$ of schemes, we get an induced functor
$$
f':V_1'\lra V_2',\quad (U_1,h_1, \shF_1,\psi_1)\longmapsto(U_1,f\circ h_1, \shF_1,\psi_1),
$$
and the formation  $f\mapsto f'$ is covariant functorial, in the strict sense.
From this one obtains  a pair of adjoint functors
\begin{equation}
\label{continuous map topoi}
f'_*:\Sh(V_1')\lra\Sh(V'_2)\quadand f'^{-1}:\Sh(V'_2)\lra \Sh(V'_1),
\end{equation}
which define a \emph{continuous map  of topoi} $\Sh(V_1')\ra\Sh(V'_2)$, in other words, 
$f'^{-1}$ commutes with finite inverse limits (\cite{SGA 4a}, Expos\'e IV, Definition 3.1).
To be explicit, the inverse image is given by
\begin{equation}
\label{effect inverse image} 
\Gamma((U,h,\shF,\psi),f'^{-1}F_2) = \Gamma((U,f\circ h,\shF,\psi), F_2),
\end{equation}
with $(U,h,\shF,\psi)\in V'_1$, whereas the direct image  takes the form 
\begin{equation}
\label{effect direct image} 
\Gamma((U,h,\shF,\psi),f'_*F_1) = \Gamma\left((U\times_{V_2}V_1,\pr_2,\pr_1^*(\shF),\pr_1^*(\psi)),F_1\right),
\end{equation}
now with $(U,h,\shF,\psi)\in V'_2$. It is straightforward to determine the effect on  structure sheaves and tautological sheaves:

\begin{proposition}
\mylabel{effect on sheaves}
In the above setting we have identifications 
$$
f'^{-1}(\O_{V'_2})=\O_{V'_1}\quadand 
f'^{-1}(\shF^\taut_{V'_2})=\shF^\taut_{V'_1} \quadand 
f'_*(\shF_{V'_1}^\taut) = \shF^\taut_{V'_2}\otimes f'_*(\O_{V'_1}).
$$
\end{proposition}

\proof
Using    \eqref{effect inverse image} with $\tilde{U}=(U,h,\shF,\psi)$ from $V'_1$ we immediately get
$$
\Gamma(\tilde{U},f^{-1}(\O_{V'_2})) = \Gamma((U,f\circ h,\shF,\psi), \O_{V'_2}) = \Gamma(U,\O_U)= \Gamma(\tilde{U}, \O_{V'_1}),
$$
and likewise for the tautological sheaves. 
Suppose now that $\tilde{U}=(U,h,\shF,\psi)$ is from $V'_2$. Consider the fiber product  $U_1=U\times_{V_2}V_1$ with respect
to $f:V_1\ra V_2$, and write $h_1:U_1\ra V_1$ for the second projection.
Using \eqref{effect direct image}   we obtain
$$
\Gamma(\tilde{U},f'_*(\shF_{V'_1}^\taut))=\Gamma( (U_1,h_1,\pr_1^*(\shF),\pr_1^*(\psi)), \shF_{V'_1}^\taut) =\Gamma(U_1, \pr_1^*(\shF) ).
$$
The latter equals $\Gamma(U,\shF\otimes\pr_{1*}(\O_{U_1}))$, as one sees by 
applying the projection formula for $\pr_1:U_1\ra U$ and the locally free sheaf $\shF$. By definition,
this coincides with the group of local sections of $\shF^\taut_{V'_2}\otimes f'_*(\O_{V'_1})$ over $\tilde{U}$.
\qed

\medskip
For each   $\tilde{U}=(U,h, \shF,\psi)$,
we get an inclusion $\GG_m\subset\uAut_{\tilde{U}/(U,h)}$ by sending an invertible scalar $\lambda$
to the automorphism $(f,\Psi)$ with $f=\id_U$ and $\Psi=\lambda\cdot\id_{\shF}$.
One easily checks that  this gives $V'\ra(\Aff/V)$ the structure of a $\GG_m$-gerbe.
In particular,   the two projections $V'\ra (\Aff/V)$ and $B(\GG_{m,V})\ra (\Aff/V)$
are \'etale locally equivalent. It follows that the fibered groupoid $V'$ is an Artin stack (\cite{Laumon; Moret-Bailly 2000}, Example 4.6.1).

Moreover, given an object $\tilde{U}=(U,h, \shF,\psi)$ from $V'$  as well as an open set   $U_0\subset U$ we get a new
object $\tilde{U}_0=(U_0,h|U_0,\shF|U_0,\psi|U_0)$ by restricting the additional data. In turn,
for each sheaf $F'$ on the stack $V'$ we get sheaves  $F_{\tilde{U}}$ on the schemes $U$.
An $\O_{V'}$-module $\shE'$ is called \emph{quasicoherent} if    the $\shE_{\tilde{U}}$ are quasicoherent sheaves on $U$  in the usual sense, for all objects $\tilde{U}\in V'$.
Particular examples are the locally free sheaves of finite rank.

For each quasicoherent sheaf $\shE'$ on $V'$, and each object $\tilde{U}=(U,h, \shF,\psi)$, we   get
an action  of the group $\Gamma(U,\GG_m)$ on  the quasicoherent sheaf $\shE'_{\tilde{U}}$ on $U$. 
By compatibility with restrictions, this  yields a \emph{linearization}
with respect to the group scheme $\GG_{m,U}$, and thus a weight decomposition
$\shE'_{\tilde{U}} =\bigoplus_{w\in\ZZ} \shE'_{\tilde{U},w}$,
as explained in \cite{SGA 3a}, \'Expose I, Proposition 4.7.3.   
This defines the \emph{weight decomposition} $\shE'=\bigoplus_{w\in\ZZ} \shE'_w$ for the sheaf on the stack $V'$.
If for a given $w\in \ZZ$ the inclusion $\shE'_w\subset\shE'$ is an equality,  one
 says that   $\shE'$ is \emph{pure of weight  $w$}.

\begin{proposition}
\mylabel{weights for sheaves}
The structure sheaf $\O_{V'}$ and the tautological sheaf $\shF^\taut_{V'}$ are pure, of respective weights $w=0$ and $w=1$.
\end{proposition}

\proof
For each object $\tilde{U}=(U,h, \shF,\psi)$ and each   $\lambda\in\Gamma(U,\GG_m)$,
the action on $\Gamma(\tilde{U},\shF_{V'}^\taut)$  is multiplication by $\lambda^w$ with $w=1$,
which one sees by making the restriction maps in \eqref{sheaves on gerbe} explicit.
Similarly, one sees that the action 
on $\Gamma(\tilde{U},\O_{V'})$ is trivial, that is, has weight  $w=0$.
\qed

\medskip
Write $\QCoh(V')$ for the    abelian category  of quasicoherent sheaves $\shE'$ and  linear maps.
Let $f:V_1\ra V_2$ be a  morphism. Note that the    linear pullback 
$$
f'^*(\shE'_2)=f'^{-1}(\shE'_2)\otimes_{f'^{-1}(\O_{V'_2})} \O_{V'_1}
$$
can be identified with the set-theoretic pullback  $f'^{-1}(\shE'_2)$, in light of Proposition \ref{effect on sheaves}.

\begin{proposition}
\mylabel{weights respected}
For each morphism $f:V_1\ra V_2$, the adjoint functors \eqref{continuous map topoi} respect   quasicoherence, 
and also   weight decompositions.
\end{proposition}

\proof
We start with the preimage functor. With $\tilde{U}=(U,h,\shF,\psi)$ from $V'_1$ we see from \eqref{effect inverse image} that 
$$
\Gamma(\tilde{U},f^{-1}\shE_2) = \Gamma((U,f\circ h,\shF,\psi),\shE_2).
$$
It follows that $(f^{-1}\shE_2)_{\tilde{U}} =(\shE_2)_{(U,f\circ h,\shF,\psi)}$ as sheaves on $U$, so quasicoherence is preserved. Moreover, if $\shE_2$ is pure of weight $w$,
the action of $\lambda\in\Gamma(U,\GG_m)$ on the group of local sections is via $\lambda^w$. By the above equality,
this also holds for $f^{-1}(\shE_2)$, which is therefore pure of weight $w$. Being left adjoint, the preimage functor respects
all direct limits, and thus direct sum decompositions. Summing up, the weight decompositions are preserved.

We   come to the direct image functor. Now let $\tilde{U}=(U,h,\shF,\psi)$ be an object from $V'_2$, and let $U_1=U\times_{V_2}V_1$ be the
base-change with respect to $f:V_1\ra V_2$. According to \eqref{effect direct image} we have
$$
\Gamma(\tilde{U},f_*\shE_1)= \Gamma((U_1,\pr_2,\pr_1^*(\shF),\pr_1^*(\psi)),\shE_1).
$$
It follows that $(f_*\shE_1)_{\tilde{U}}= \pr_{1,*} ((\shE_1)_{(U_1,\pr_2,\pr_1^*(\shF),\pr_1^*(\psi))})$ as sheaves on $U$.
Since the schemes $U,V_1,V_2$ are affine, the projection  $\pr_1:U_1=U\times_{V_2}V_1\ra U$  is affine, so
the direct image functor $\pr_{1,*}$ preserves quasicoherence. Thus the sheaf $(f_*\shE_1)_{\tilde{U}}$  on $U$ is quasicoherent.
Moreover, for quasicoherent sheaves the  functor $\pr_{1,*}$ preserves direct sum decompositions.
If $\shE_1$ is pure of weight $w$, we argue as in the preceding paragraph to see that $f_*(\shE_1)$ is pure of the same weight.
In turn, the weight decompositions are preserved.
\qed

%===========================================================
\section{The stack of twisted Schur representations}
\mylabel{Stacks schur}

Let $R$ be a ground ring, $\Lambda$ be an associative algebra, and $n\geq 1$ be some integer.
We call a pair $(E,\rho)$, where $E$ is a locally free $R$-module of rank $n$, and $\rho:\Lambda\ra\End(E)$ is a homomorphism
of associative algebras, a \emph{linear representation  of degree $n$}. For some suitable Zariski open covering $A_i=R[1/f_i]$ 
and   choices of bases for $E\otimes_RA_i$, 
this yields a collection of  \emph{matrix representations} $\rho_i:\Lambda\otimes_R A_i\ra \Mat_n(R)\otimes_R A_i$, not necessarily
compatible.
The goal of this section is to  clarify these seemingly innocuous  facts by using stacks, and 
generalize it from matrix algebras to Azumaya algebras. Furthermore, we relate it 
to our results on modifying moduli problems and  spaces of representations  in  
Sections \ref{Points geometric origin} and \ref{Spaces representations}.
Indeed,   Azumaya algebras are needed to perform and explain such modifications.

Throughout,  we fix an Azumaya algebra $\Lambda^\azu$  of   degree $n\geq 1$,
with $\Lambda^\azu=\Mat_n(R)$ as important special case.  
To apply the results of Section \ref{Spaces representations},
we also assume that the given associative algebra $\Lambda$ is finitely presented.
Consider the resulting functor of representations 
$$
F:(\Aff/R)\lra(\Set),\quad A\longmapsto\Hom_{\text{$A$-Alg}}(\Lambda\otimes_RA,\Lambda^\azu\otimes_RA).
$$
From Theorem \ref{space  schur representations} we see that the subfunctor of Schur representations is representable by 
a quasiaffine scheme of finite presentation,  which we denote by $X=X^{\Lambda^\azu}_{\Lambda/R}$.
 
Recall that we have a short exact sequence 
$$
1\lra \GG_m\lra U_{\Lambda^\azu/R}\stackrel{\conj}{\lra} \Aut_{\Lambda^\azu/R}\lra 1,
$$
where the terms on the right are the respective group schemes of units and automorphisms for $\Lambda^\azu$.
According to Proposition \ref{schur subfunctor}, the canonical action of $\Aut_{\Lambda^\azu/R}$ on the sheaf $F$ stabilizes the subsheaf  
of Schur representations, and the induced action on the  quasiaffine scheme $X$ is free. 
To conform with the notation in Section \ref{Quotient stacks}, we  now set $G=U_{\Lambda^\azu/R}$ and $H=\Aut_{\Lambda^\azu/R}$.
The quotient $Q=X/H^\op$ is representable by an algebraic space, and its formation commutes with 
base-change along arbitrary $\tilde{Q}\ra Q$, see for example \cite{Laurent; Schroeer 2024}, Lemma 1.1.
Note that such quotients may easily be \emph{non-separated}, the simplest example   stemming
from the  action $\lambda\cdot (x_1,x_2)=(\lambda x_1,\lambda^{-1}x_2)$ of the multiplicative group on the pointed affine plane, whose
quotient is the affine line with double origin. 

Also note that such quotients are usually \emph{non-schematic}. Examples are abundant, but at the same time far from obvious:
Working over a ground field $k=k^\alg$, we start with some
non-schematic   surface $Q$ that is proper, integral and normal (\cite{Artin 1973}, Example 4.4).
This has the resolution property (\cite{Mathur; Schroeer 2021}, Theorem 6.8; the arguments in  \cite{Schroeer; Vezzosi 2004}, which 
hold true for algebraic spaces, already suffice). 
According to Totaro's result (\cite{Totaro 2004}, Theorem 1.1),
there is a locally free sheaf $\shE$  such that the corresponding principal bundle $P\ra Q$
with respect to $G=\GL_n$, $n=\rank(\shE)$ has quasiaffine total space. The arguments in 
loc.\ cit.\ (below the proof of Corollary 5.2) reveal that we may replace the sheaf by  $\shE\oplus\shE^\vee$, and
thus may assume that the principal bundle admits a reduction of structure $P_0\subset P$ to $G_0=\SL_n$, which contains
$\mu_n$ as a finite normal subgroup scheme.
Setting $X=P_0/\mu_n$, we get a principal bundle  $X\ra Q$  with respect to $H=\PGL_n$.
By construction, $P_0$ is quasiaffine, $X$ is normal and schematic, and the quotient map $P_0\ra X$  is finite and flat.
Viewing $\O_X$ as the norm of $\O_{P_0}$, we see that $\O_X$ is ample (\cite{EGA II}, Proposition 6.6.1).
Summing up, $X$ is quasiaffine with a free action of $H=\PGL_n$  such that the quotient $Q$
is non-schematic.

There are also examples taken from the realm of moduli spaces:
Wei\ss mann and Zheng showed that for every smooth proper curve of genus  $g\geq 4$ the coarse moduli space
of simple  sheaves that are locally free of rank $n\geq 2$ is a non-schematic algebraic space (\cite{Weissmann; Zhang 2023}, Corollary 2.8).
 
Back to our $X=X^{\Lambda^\azu}_{\Lambda/R}$.
The algebraic space $Q=X/H^\op$ can be regarded as a ``moduli space'' $M=M^{\Lambda^\azu}_{\Lambda/R}$ of Schur representations.
However, we regard such locutions
as  ``dangerous'', and  seek to define the ``true'' \emph{moduli stack} $\shM=\shM^{\Lambda^\azu}_{\Lambda/R}$ of Schur representations,
in terms of ``concrete'' representation-theoretic data.
We then construct a \emph{comparison functor} $\Phi:\shM\ra [X/G^\op/Q]$ to the quotient stack, which is an equivalence relating 
representation-theoretic and  algebro-geometric data. The quotient space  $Q=X/H^\op$ turns out to be the  \emph{coarse moduli space}.

Throughout, we find it  psychologically helpful to regard  schemes, which by definition are certain ringed spaces $V=(|V|,\O_V)$,
not only as    functors $(\Aff/ R)\ra(\Set)$ via the Yoneda embedding,
but actually    as   \emph{fibered groupoids} $(\Aff/V)\ra(\Aff/  R)$ via the comma construction.
Recall that   $V'\ra (\Aff/V)$ denotes the splitting gerbe for the Azumaya algebra $\Lambda^\azu\otimes_R\O_V$.
We come to a central notion:

\begin{definition}
\mylabel{twisted representation}
A \emph{twisted representation} of the associative algebra $\Lambda$ in the Azumaya algebra $\Lambda^\azu$   over a scheme $V$ 
is a pair $(\shE',\rho')$ where 
$\shE'$ is  a locally free sheaf   on  the $\GG_{m,V}$-gerbe $V'$  that has rank $n=\deg(\Lambda^\azu)$ and is pure of weight $w=1$, 
and $\rho':\Lambda\otimes_R\O_{V'}\ra \uEnd(\shE')$ is a homomorphism of $\O_{V'}$-algebras. 
\end{definition}

If for all objects $\tilde{U}=(U,h,\shF,\psi)$ from $V'$ the   $\rho'_{\tilde{U}}:\Lambda\otimes_R\O_U\ra\uEnd(\shE'_{\tilde{U}})$ 
are Schur representations,
we say that  $(\shE',\rho')$ is a \emph{twisted Schur representation}.

Let us emphasize that we regard ``twisted'' representations as ``true'' representations, albeit defined on a stack rather than  a scheme.
In other words, we favor the approach of de Jong  and Lieblich  (\cite{de Jong 2006}, \cite{Lieblich 2004}, \cite{Lieblich 2007})   
over C\u{a}ld\u{a}raru's point of view (\cite{Caldararu 2000},  \cite{Caldararu 2002}).
Consider now  the  category  
$$
\shM=\shM^{\Lambda^\azu}_{\Lambda/R}= \{(V,\shE',\rho')\mid \text{$\rho':\Lambda\otimes_R\O_{V'}\ra \uEnd(\shE')$}\}
$$
whose objects  are triples  $(V,\shE',\rho')$ where $V$ is an affine scheme 
and $(\shE',\rho')$ is a twisted Schur representation over the scheme $V$, in other words, a weight-one Schur representation over the stack $V'$.
Morphisms $(V_1,\shE'_1,\rho'_1)\ra (V_2,\shE'_2,\rho'_2)$   are pairs $(f,\Psi')$ where $f:V_1\ra V_2$ is a morphism of schemes,
and $\Psi':f'^{-1}(\shE'_2)\ra\shE'_1$ is an isomorphism of modules over $\O_{V'_1}=f^{-1}(\O_{V'_2})$ making the diagram
$$
\begin{tikzcd} 
			& \Lambda\otimes_R\O_{V_1'}\ar[dl,"f'^{-1}(\rho'_2)"']\ar[dr,"\rho'_1"]\\
\uEnd(f'^{-1}(\shE'_2))  \ar[rr,"\Psi'^{\vee-1}\otimes \Psi'"'] 	& 				&  \uEnd(\shE'_1)
\end{tikzcd}
$$
commutative.  

\begin{proposition}
\mylabel{shm fibered groupoid}
The forgetful functor $\shM\ra (\Aff/R)$ given by $(V,\shE',\rho')\mapsto V$ endows $\shM$ with the structure of a category fibered in groupoids.
Moreover, the stack axioms hold with respect to the fppf topology.
\end{proposition}

\proof
By definition, in morphisms $(f,\Psi'):(V_1,\shE'_1,\rho'_1)\ra (V_2,\shE'_2,\rho'_2)$ the linear map $\Psi'$ is an isomorphism.
From this one immediately infers that all morphisms in $\shM$ are cartesian. In particular, all fiber categories are groupoids
(\cite{SGA 1}, Expos\'e VI, Remark after Definition 6.1; but note that the very existence of cartesian maps was overlooked there).
We next check that  arrows lift with respect to  the forgetful functor,  with prescribed target. This is immediate:
Suppose that $f:V_1\ra V_2$ is a morphism of affine schemes,
and $(V_2,\shE'_2,\rho'_2)$ is an object   over $V_2$.
Setting $\shE'_1=f^{-1}(\shE'_2)$ and $\rho'_1=f^{-1}(\rho'_2)$ and $\Psi'=\id_{\shE'_1}$, we get the desired morphism $(f,\Psi')$ over $f$.

It remains to verify that $\shM$ is a stack. We first check that the automorphism presheaves satisfy the sheaf axioms.
Let $A$ be an $R$-algebra,   $A\subset A_0$ be an fppf extension, and write $A_1=A_0\otimes_AA_0$.
Set $V_i=\Spec(A_i)$ and fix an object $(V,\rho',\shE')$ over $V=\Spec(A)$.
We have to check that $\Aut ((V,\rho',\shE')/V)$ is an equalizer of the diagram
$$
\begin{tikzcd} 
\Aut ((V_0,h_0^{-1}(\rho'),h_0^{-1}(\shE'))/V_0) \ar[r, shift left=.5ex] \ar[r, shift right = .5ex]	& \Aut ((V_1,h_1^{-1}(\rho'),h_1^{-1}(\shE'))/V_1)
\end{tikzcd}
$$
where $h_i:V_i\to V$ denote the canonical morphisms. First suppose $( \id_V ,\Psi')$ is an automorphism over $V$ that becomes the 
identity over $V_0$.  Hence  $\Psi'$ is an automorphism of $\shE'$ compatible with  $\rho'$.
Consider the quasicoherent sheaves $\shE'_{\tilde{U}}$ on the affine scheme $U$, which are indexed by 
the objects $\tilde{U}=(U,h,\shF,\psi)$ of the splitting gerbe $V'$.
Our $\Psi'$ is determined by the compatible collection of automorphisms $\Psi'_{\tilde{U}}:\shE'_{\tilde{U}}\ra\shE'_{\tilde{U}}$.
The latter become identities on $U\times_{V_0}V_1$,  so by  \cite{SGA 1}, Expos\'e VIII, Theorem 1.1
the  $ \Psi'_{\tilde{U}}$ must be identities.
Likewise, one checks that every automorphism over $V_0$ whose two inverse images  on $V_1$ coincide
comes from an automorphism over $V$.

 Finally we need  to check that every descent datum is effective.
Using the above notation, we now  fix an object $\zeta_0=(V_0,\shE'_0,\rho'_0)$ over $V_0=\Spec(A_0)$. 
Consider the two inclusions
$A_0\rightrightarrows A_1$ given by $x\mapsto 1\otimes x$ and $x\mapsto x\otimes 1$. As customary, we 
write $\zeta_0\otimes_{A_0}A_1$ and $A_1\otimes_{A_0}\zeta_0$ for the resulting pullbacks, say    
formed as  in the first paragraph above. Suppose now that we have an isomorphism
$$
(f_1,\Psi'_1):\zeta_0\otimes_{A_0}A_1 \ra A_1\otimes_{A_0}\zeta_0
$$
that satisfies the cocycle condition over $V_2=\Spec(A_2)$, with $A_2=A_0\otimes_A A_0\otimes_AA_0$.
Consider the   collection of locally free sheaves $\shE'_{\tilde{U}_0}$ arising from $\shE'_0$, indexed by the objects
$\tilde{U}_0=(U_0,h_0,\shF_0,\psi_0)$ of the splitting gerbe $V'_0$. These live  on the affine schemes $U_0$, and are endowed
with   comparison morphisms with respect to the morphisms in $V'_0$.
Applying \cite{SGA 1}, Expos\'e VIII, Theorem 1.1 for each affine scheme $U$ appearing in the objects  $\tilde{U}=(U,h,\shF,\psi)$ of $V'$, 
we get a   collection $\shE'_{\tilde{U}}$ together with comparison maps.  
The desired sheaf on $V'$  is now defined via $\Gamma(\tilde{U},\shE')=\Gamma(U,\shE'_{\tilde{U}})$, where the restriction
maps stem from the comparison maps.
\qed

\medskip
We call $\shM=\shM^{\Lambda^\azu}_{\Lambda/R}$ the \emph{stack of twisted   Schur representations} of the 
associative algebra $\Lambda$ in the Azumaya algebra $\Lambda^\azu$.
The task now is to construct a comparison  functor 
$\Phi:\shM \ra [X/G^\op/Q] $
that relates representation-theoretic   with algebro-geometric data. Recall that the quotient stack comprises
tuples $(U,g,P,f)$  where $U$ is an affine scheme, $g:U\ra Q$ is a morphism of algebraic spaces,
$P$ is a $G_U$-torsor, and  $f:P\ra X_U$ is a $G_U$-equivariant morphism that   induces $g$ on quotients.

Let $(V,\shE',\rho')$ be an object from   $\shM$. On  the gerbe of splittings $V'$ we then have two locally free sheaves
$\shE'$ and $ \tautsheaf$, both having rank $n=\deg( \Lambda^\text{azu})$ and weight $w=1$.
Composition endows the sheaf  $\uHom(\shE', \tautsheaf )$ with a module structure over $\uEnd(  \tautsheaf )$,
and the latter  comes with an isomorphism $ \tautsplit :\uEnd( \tautsheaf )\ra \Lambda^\text{azu}\otimes_R\O_{V'}$.
These sheaves are locally free of rank $n^2$ and weight $w=0$, and thus correspond to locally free sheaves of
likewise rank on the affine scheme $V$.   Let us write
$$
P_{V,\shE'}\lra V
$$
for the resulting   torsor with respect to the group scheme $G_V=U_{\Lambda^\azu/R}|V$; 
it corresponds to  the locally free module of rank one over the sheaf of Azumaya algebras $\Lambda^\azu\otimes_R\O_{V'}$ 
given by $\uHom_{\O_{V'}}(\shE', \tautsheaf )$.
Note that the base and the total space of the torsor are affine.
One easily checks that the formation of $P_{V,\shE'}$  is functorial in $(V,\shE',\rho')$.

Our next task is to construct a  morphism of schemes   $f= f_{V,\shE',\rho'}$
from   $P=P_{V,\shE'}$ to the quasiaffine scheme $X=X^{\Lambda^\azu}_{\Lambda/R}$ 
representing the functor $F^0$ of Schur representations of $\Lambda$ in $\Lambda^\azu$.
 Note that the former is affine, whereas the latter is quasiaffine.
We shall specify  the morphism as a  functor $(\Aff/P)\ra (\Aff/X)$.
Let $U$ be an affine scheme, together with a morphism $U\ra P$. By the universal property of fiber products, this  can be seen
as a morphism $f:U\ra V$, together with a section of the induced $G_U$-torsor $P_U=P\times_VU$.
By the very definition of the torsor, the 
latter is nothing but an isomorphism $\psi':\shE'|U'\ra  \tautsheaf[U']$ of locally free sheaves on the gerbe of splittings $U'$.
Now the given linear representation $\rho':\Lambda\otimes_R\O_{V'}\ra \uEnd(\shE')$ enters the picture:
The composite map 
$$
\Lambda\otimes_R\O_{U'}\stackrel{f'^{-1}(\rho')}{\lra} \uEnd(f'^{-1}\shE') 
\stackrel{\psi'^{\vee-1}\otimes\psi'}{\lra} \uEnd(  \tautsheaf[U'])\stackrel{\psi^\taut_{U'}}{\lra}\Lambda^\azu\otimes_R\O_{U'}
$$
is a Schur representation of $\Lambda$ in $\Lambda^\azu$ over the splitting gerbe $U'$.
It corresponds to a Schur representation over the scheme $U$, because the 
involved sheaves have weight $w=0$.  Since $X=X^{\Lambda^\azu}_{\Lambda/R}$ represents the functor $F^0$, this can be seen as a morphism $ U\ra X$.
One easily checks that the formation is functorial, and thus defines the desired morphism  
$$
f_{V,\shE',\rho'}:P_{V,\shE'}\lra X.
$$
Recall that the group scheme $G=U_{\Lambda^\azu/R}$ acts   on $X$  via   conjugation, 
and on the $G_V$-torsor $P_{V,\shE'}$ in an obvious way. The following technical fact is a crucial   observation: 

\begin{lemma}
\mylabel{equivariant}
The above morphism $f_{V,\shE',\rho'}:P_{V,\shE'}\ra X$ is equivariant with respect to the $G$-actions.
\end{lemma} 
 
\proof
We have to check that the map $P_{V,\shE'}(A)\ra X(A)$ is equivariant for the abstract groups $G(A)$,
for each $R$-algebra $A$. The morphism $P_{V,\shE'}\ra V$ endows $A$ with an algebra structure over $A_0=\Gamma(V,\O_V)$.
Since our constructions commute with base-change, it suffices to treat the case $A=A_0=R$. 

Fix some $\sigma\in G(R)=(\Lambda^\azu)^\times$.  Recall that its effect on the set $P_{V,\shE'}(R)=\Hom(\shE',  \tautsheaf )$ is via composition
with the  automorphism  $\eta'=(  \tautsplit )^{-1}(\sigma)$  from 
$\End( \tautsheaf )$. By definition, our map $f_{V,\shE',\rho'}$ sends the element of $P_{V,\shE'}(R)$ 
corresponding to a linear map $\psi':\shE'\ra  \tautsheaf $
to the element  $X(R)$ stemming from the composite map
$$
\Lambda\otimes_R{\O_{V'}}\stackrel{\rho'}{\lra}\uEnd(\shE')\stackrel{\psi'^{\vee  -1}\otimes\psi'}{\lra} 
\uEnd(  \tautsheaf )\stackrel{\psi^\taut_{V'}}{\lra} \Lambda^\azu\otimes_R\O_{V'}.
$$
So $f_{V,\shE',\rho'}(\sigma\cdot\psi')$ is given by a likewise  composition, formed with 
$$
(\eta'\circ\psi')^{\vee -1}\otimes (\eta'\circ\psi') = (\eta'^{\vee -1}\otimes\eta') \circ (\psi'^{\vee -1}\otimes\psi')
$$
instead of $ \psi'^{\vee -1}\otimes\psi'$. The task is to verify   that composition of   endomorphisms of $ \tautsheaf$ 
with   $\eta'^{\vee -1}\otimes\eta'$ is the same as conjugation with $\eta'$.
This has to be  checked over the objects   $(U,h,\shF,\varphi)$ of the splitting gerbe $V'$,  thus becomes 
a statement about sheaves on the affine schemes $U$, hence, a problem in  commutative algebra:

Suppose $E$ is a locally free $R$-module of rank $n$. Recall that the canonical map 
\begin{equation}
\label{cartan identification}
E^\vee\otimes E\lra \End(E),\quad f\otimes a\longmapsto (x\mapsto f(x)\cdot a).
\end{equation}
is bijective, and that under this identification     composition of endomorphisms corresponds to the pairing of tensors  
$(f\otimes a)\otimes (g\otimes b)\mapsto f(b)\cdot g\otimes a$, see \cite{A 1-3}, Chapter II, \S 4, No.\ 2.
Each $h\in \GL(E)$ yields  bijective linear maps
$$
E^\vee\otimes E\lra E^\vee\otimes E\quadand \End(E)\lra \End(E)
$$
given by $h^{\vee-1}\otimes h$ and $g\mapsto h\circ g\circ h^{-1}$, respectively.
We have  to verify that these bijections coincide   under   \eqref{cartan identification}.
The problem is local, so we may assume that there is a basis $e_1,\ldots,e_n\in E$,
with ensuing identification
$\End(E)=\Mat_n(R)$, standard basis  $E_{ij}\in\Mat_n(R)$, and dual basis $e_1^\vee,\ldots,e_n^\vee\in E^\vee$.
Write  $(\alpha_{ij})$ and $(\beta_{ij})$ for  the matrices for $h$ and $h^{-1}$, respectively. 
Note that   the tensor $ e_s^\vee\otimes e_r$ corresponds to the matrix $E_{rs}$, and that $h^{\vee-1}=h^{-1\vee}$ has matrix
$(\beta_{ji})$.
One computes
$$
(h^{\vee-1}\otimes h)(e_s^\vee\otimes e_r) = (\sum_j \beta_{sj} e^\vee_j)\otimes(\sum_i\alpha_{ir} e_i)=\sum_{i,j}(\beta_{sj}\alpha_{ir}\cdot e_j^\vee\otimes e_i).
$$
Using Kronecker deltas, we write $E_{rs}=(\delta_{ir} \delta_{sj})$ and see that the $(i,j)$-entry of
the triple matrix product $(\alpha_{ik})(\delta_{ir} \delta_{sj}) (\beta_{kj})$ is given by 
$\sum_{k,l}\alpha_{ik}\cdot\delta_{kr}\delta_{sl}\cdot\beta_{lj} = \alpha_{ir}\beta_{sj}$. Thus
$h\circ E_{rs}\circ h^{-1} = \sum_{i,j} \alpha_{ir}\beta_{sj}\cdot E_{ij}$,
which indeed corresponds to $(h^{\vee-1}\otimes h)(e_s^\vee\otimes e_r)$.
\qed

\medskip
Since the  $G_V$-action on $P_{V,\shE'}$ is free, our equivariant $f_{V,\shE',\rho'}$ induces a morphism 
$$
g_{V,\shE',\rho'}:V=P/   G_V^\op  \lra X/G^\op = X/H^\op =Q.
$$
Summing up, we have attached to every  object $(V,\shE',\rho')\in\shM$ from the stack of   twisted Schur representations an object 
$$
(V,g_{V,\shE',\rho'},P_{V,\shE'},f_{V,\shE',\rho'})\in[X/G^\op/Q]
$$
from the quotient stack. One easily checks that the formation is functorial. This defines  the desired \emph{comparison functor}
$$
\Phi:\shM=\shM^{\Lambda^\azu}_{\Lambda/R}\lra [X/G^\op/Q].
$$
Obviously, this is compatible with the forgetful functor to $(\Aff/R)$.
One main result of this paper is that this relation between representation-theoretic   and algebro-geometric data 
is essentially an identification:

\begin{theorem}
\mylabel{comparison functor equivalence}
The comparison functor $\Phi:\shM\ra [X/G^\op/Q]$ is an equivalence of categories.
\end{theorem} 

\proof
We have to check that $\Phi$ is fully faithful and essentially surjective. 
This relies to a large degree on  internal properties of
the category $\shM$.
As a preparation, note that we also have a functor
\begin{equation}
\label{category over Q}
\shM\lra (\Aff/Q),\quad (V,\shE',\rho')\longmapsto (V,g_{V,\shE',\rho'}).
\end{equation}
So both $\shM$ and $[X/G^\op/Q]$ are categories over $(\Aff/Q)$, and  $\Phi$ is a functor relative to this.
We  now proceed in three steps.
 
\medskip
{\bf Step 1:} 
\emph{The above functor \eqref{category over Q} endows $\shM$ with the structure of a category  fibered in groupoids.
Moreover, the stack axioms hold with respect to the fppf topology on   $(\Aff/Q)$.}
The arguments are exactly  as in the proof for Proposition \ref{shm fibered groupoid}.

\medskip
{\bf Step 2:} 
\emph{The comparison functor $\Phi$ is fully faithful.}
In light of \cite{SGA 1}, Expos\'e VI,  Proposition 6.10 and step 1,  it suffices to verify this
for the fiber categories over $(\Aff/Q)$.
Fix an affine scheme $V$ and a  morphism of algebraic spaces $g:V\ra Q$. Consider
an object  $(V,\shE',\rho')$  in the stack of Schur representations $\shM$ over $(V,g)$,
in other words $g=g_{V,\shE',\rho'}$. Setting  $P=P_{V,\shE'}$ and $f=f_{V,\shE',\rho'}$ we get an object  $(V,g,P,f)$
in the quotient stack  $[X/G^\op/Q]$ over $(V,g)$.  
According to Corollary \ref{abelian gerbe}, the structure map $R\ra\Lambda^\azu$ induces a bijection
$\Gamma(V,\O_V^\times)\ra \Aut_{(V,g)}(V,g,P,f)$.
Scalar multiplication of the structure sheaf on $\shE'$ gives an injective  
group homomorphism $\Gamma(V,\O_V^\times)\ra \Aut_V (V,\shE',\rho')$. 
%By Proposition \ref{shm fibered groupoid} the   endomorphism monoid coincides with the automorphism group. 
Since the representation $\rho'$
is Schur, the injection $\Gamma(V,\O_V^\times)\ra \Aut_V (V,\shE',\rho')$ is actually bijective.
Consider the diagram
$$
\begin{tikzcd} [row sep=tiny]
			%& \Aut_{(V,g)}(V,g,P,f)\ar[dd,"\Phi"]\\
			&  \Aut_V (V,\shE',\rho') \ar[dd,"\Phi"]\\
\Gamma(V,\O_V^\times)\ar[ur]\ar[dr]\\
			& \Aut_{(V,g)}(V,g,P,f),
			%& \Aut_V (V,\shE',\rho'),
\end{tikzcd}
$$
where the diagonal arrows are as specified above.
Using that the $G_V$-torsor $P$ stems from the locally free sheaf $\uHom(\shE', \tautsheaf )$,
we infer that the diagram is commutative.
It follows that the vertical arrow is bijective, and moreover that every automorphism
of $(V,\shE',\rho')$ over $V$ is actually over $(V,g)$. 
So $\Phi$ yields bijections on automorphism groups in fiber categories. Since the fibers are groupoids,
$\Phi$ is fully faithful on the fibers.
 
\medskip
{\bf Step 3:} 
\emph{The   functor $\Phi$ is essentially surjective.} 
In both $\shM$ and $[X/G^\op/Q]$, viewed as fibered categories over $(\Aff/Q)$,  every descent datum  is effective.
It thus suffices to check that every descent datum of the latter arises from a descent datum on the former.
But this is immediate from step 2.
\qed

\medskip
Recall that for any fibered groupoid $\shG$ over $(\Aff/R)$, a morphism $\shG\ra Z$ to an algebraic space $Z$ that is universal for morphisms 
into algebraic spaces is called   \emph{coarse moduli space}.

\begin{corollary}
\mylabel{twisted representation artin}
The fibered groupoid $\shM\ra(\Aff/R)$   is an Artin stack, and  the morphism $\Phi^\text{\rm crs}:\shM\ra Q=X/H^\op$ is the  coarse moduli space.
Moreover, $\shM$ carries the structure of a $\GG_m$-gerbe over $(\Aff/Q)$.
\end{corollary}

\proof
It suffices to check these statements for the quotient stack $[X/G^\op/Q]$, which by the theorem is equivalent to $\shM$, 
as categories over $(\Aff/Q)$.
We saw in  Section \ref{Quotient stacks} that the quotient stack is fibered in groupoids, 
satisfies the stack axioms, and is a $\GG_m$-gerbe over $(\Aff/Q)$.
According to \cite{Laumon; Moret-Bailly 2000}, Example 4.6.1 it is indeed an Artin stack. 
It remains to check the statement on the coarse moduli space. Obviously, the transformation
$$
[X/G^\op/Q]\lra [X/H^\op/Q]=Q,\quad (U,g,P,f)\longmapsto (U,g,\bar{P},\bar{f}),
$$
where $\bar{P}= H\wedge^GP  = P/\GG_m$ is the induced $H$-torsor, and $\bar{f}:\bar{P}\ra X_U$ is the induced   $H_U$-equivariant map, 
yields  the sheafification of the functor of isomorphism
classes for the $\GG_m$-gerbe. It follows that $\shM\ra Q$ is the universal morphism from the stack to an algebraic space.
\qed

%===========================================================
\section{Modifying moduli of representations}
\mylabel{Modifying moduli}

\newcommand{\qc}{\operatorname{qc}}
We keep the set-up of the preceding section: Let $R$ be a ground ring,  $\Lambda$ be a finitely presented  associative algebra,
and $\Lambda^\azu$ be an Azumaya algebra 
of degree $n\geq 1$, with resulting group schemes $G=U_{\Lambda^\azu/R}$ and $H=\Aut_{\Lambda^\azu/R}$  
of units and automorphisms, respectively.

In order to deal with  moduli of representations   of $\Lambda$ in $\Lambda^\azu$, 
a first attempt is to consider  the quasiaffine scheme $X=X^{\Lambda^\azu}_{\Lambda/R}$  of Schur representations.
This  ``over-represents'' our moduli problem,  because it does not take into account  the  isomorphism relation, and is therefore unsatisfactory.
One may next form the algebraic space $Q=X/H^\op$.  But this ``under-represents'' the moduli problem, since it neglects
automorphisms, and is  still unsatisfactory. The correct framework  is the quotient stack $[X/G^\op/Q]$.
This    general construction of algebraic geometry, however, obscures the representation-theoretic content of the given moduli
problem. This is rectified by Theorem \ref{comparison functor equivalence}, which   tells us that one can work with the Artin stack
$\shM=\shM^{\Lambda^\azu}_{\Lambda/R}$ of twisted Schur representations instead. 

The quasiaffine scheme $X$, the algebraic space $Q$,
and the Artin stack $\shM$ are related by  a commutative diagram of functors
\begin{equation}
\label{scheme space stack}
\begin{tikzcd} 
        	&   (\Aff/X)\ar[dl,"\can"']\ar[dr,"\can"]\\
\shM\ar[rr,"\Phi^\crs"']	&			& (\Aff/Q).
\end{tikzcd}
\end{equation} 
The diagonal arrow on the right stems from  the quotient map for $Q=X/H^\op$, 
the diagonal arrow on the left is the classifying map for the universal Schur representation $\Lambda\otimes_R\O_X\ra\Lambda^\azu\otimes_R\O_X$ 
pulled back to the splitting gerbe $X'$, and the horizontal functor is given by \eqref{category over Q}, which stems from 
the comparison morphism $\Phi:\shM\ra[X/G^\op/Q]$.

\emph{We can now easily
answer the  question  raised
at the end of Section \ref{Points geometric origin} on   modified moduli problems for the situation at hand:}
Let    $g\in Q(R)$ be an   $R$-valued point, not necessarily of geometric origin, and
consider the   $H$-torsor $g^*(X)$ and the resulting  twisted forms 
$$
\tilde{H}\quadand \tilde{X}\quadand \tilde{\Lambda}^\azu
$$
of $H$ and $X$ and $\Lambda^\azu$, respectively.  Then $\tilde{\Lambda}^\azu$ is another Azumaya algebra of degree $n$.
According to Theorem  \ref{twist geometric origin}, our $R$-valued point $g\in Q(R)$ has acquired geometric origin
with respect to $\tilde{X}/\tilde{H}^\op=Q$.    The following    unravels the 
representation-theoretic content of this statement:

\begin{theorem}
\mylabel{representation-theoretic content}
The twisted form $\tilde{X}$ is the scheme of Schur representations of $\Lambda$ in the Azumaya algebra $\tilde{\Lambda}^\azu$,
and there is a Schur representation $\tilde{\rho}:\Lambda\ra\tilde{\Lambda}^\azu$ inducing $g\in Q(R)$.
\end{theorem}

\proof
The first statement is a consequence of Proposition \ref{schur subfunctor}, and the second follows from 
Theorem  \ref{twist geometric origin}.
\qed

\medskip
Note that the above result does not rely  whatsoever  on the Artin stacks
$$
\shM=\shM^{\Lambda^\azu}_{\Lambda/R}\quadand \tilde{\shM}=\shM^{\tilde{\Lambda}^\azu}_{\Lambda/R}
$$
of twisted Schur representations. However, it raises the   question  whether the existence
of an $R$-valued object in $\shM$ already implies the existence of an $R$-valued point in $X$,
and likewise for $\tilde{\shM}$ and $\tilde{X}$.

Recall that  $S'$ denotes the gerbe of splittings for the Azumaya algebra $\Lambda^\azu$ over the base scheme $S=\Spec(R)$.
We say that a  
twisted Schur representation $\rho':\Lambda\otimes_R\O_{S'}\ra \uEnd(\shE')$ \emph{induces a given $R$-valued point $g\in Q(R)$}
if the object $(S,\shE',\rho')\in \shM(R)$   maps to $g\in Q(R)$ in the   diagram \eqref{scheme space stack}.
Note  that this generalizes our terminology from   Schur representation $\rho:\Lambda\ra\Lambda^\azu$
to the twisted case.

\begin{theorem}
\mylabel{characterization geometric origin}
Let  $g\in Q(R)$, and assume that  the non-abelian cohomology set $H^1(S,\GL_n)$ is a singleton. Then the following three conditions are equivalent:
\begin{enumerate}
\item The $R$-valued point $g\in Q(R)$ has geometric origin.
\item There is a Schur representation $\Lambda\ra\Lambda^\azu$ inducing $g$.
\item There is a twisted Schur representation $\Lambda\otimes_R\O_{S'}\ra \uEnd(\shE')$ inducing $g$.
\end{enumerate}
Moreover,   (i)$\Leftrightarrow$(ii)$\Rightarrow$(iii)  hold without the assumption on $H^1(S,\GL_n)$.
\end{theorem}

\proof
%The implications (i)$\Leftrightarrow$(ii)$\Rightarrow$(iii) are trivial indeed.
For the equivalence (i)$\Leftrightarrow$(ii) recall that $g$ is of geometric origin if and only if it is in the image of $X(R)\to Q(R)$. The implication (ii)$\Rightarrow$(iii) is trivial.
Suppose now that the set  $H^1(S,\GL_n)$ is a singleton, and that $\rho':\Lambda\otimes_R\O_{S'}\ra \uEnd(\shE')$ is a twisted representation inducing $g$.
This gives an object $(S,\shE',\rho')\in \shM$, and via the comparison functor $\Phi:\shM\ra[X/G^\op/Q]$ a $G$-torsor $P\ra S$ together with a $G$-equivariant
morphism $f:P\ra X$ making the diagram
$$
\begin{CD}
P	@>f>>	X\\
@VVV		@VVV\\
S	@>>g>	Q
\end{CD}
$$
commutative. The identifications $H^1(S,\GL_n)=H^1(S,G)$ stemming from the torsor translation maps \eqref{torsor tranlation map}
reveals that our torsor admits a section $s:S\ra P$. The composite $f\circ s:S\ra X$ corresponds to a Schur representation $\rho:\Lambda\ra\Lambda^\azu$ inducing $g$,
which therefore has geometric origin.
\qed

%===========================================================
\section{Tautological sheaves}
\mylabel{Tautological sheaves}
 
We keep the set-up of the preceding sections: $R$ is a ground ring, $\Lambda$ is a  finitely presented associative algebra, 
$\Lambda^\azu$ is an Azumaya algebra of degree $n\geq 1$, and $H=\Aut_{\Lambda^\azu/R}$.
We then have the quasiaffine scheme $X=X^{\Lambda^\azu}_{\Lambda/R}$ of Schur representations, 
and the Artin stack $\shM=\shM^{\Lambda^\azu}_{\Lambda/R}$ of twisted Schur representations, which is a $\GG_m$-gerbe
over the algebraic space $Q=X/H^\op$.
The goal of this section is to  define on the stack the  \emph{tautological tilting sheaf} $\shT_\shM$, and relate it to various
Azumaya algebras and Brauer classes.

Recall that the objects $\hat{V}=(V,\shE',\rho')$ from $\shM$ are triples $(V,\shE',\rho')$ where $V$ is an affine scheme,
$\shE'$ is a locally free sheaf of rank $n=\deg(\Lambda^\azu)$ and weight $w=1$ on the gerbe  $V'$ of splittings for $\Lambda^\azu\otimes_R\O_V$,
and $\rho':\Lambda\otimes_R\O_{V'}\ra\uEnd(\shE')$ is a Schur representation.
Also recall that the objects of $V'$ are quadruples $(U,h,\shF,\psi)$ where $U$ is an affine scheme, $h:U\ra V$ is a morphism,
$\shF$ is a locally free sheaf of rank $n$ over $U$, and $\psi:\uEnd(\shF)\ra\Lambda^\azu\otimes_R\O_U$ is an isomorphism.

We  regard $\shM$ as a site,  where the covering families $(V_\lambda,\shE'_\lambda,\rho'_\lambda)_{\lambda\in L}\ra (V,\shE',\rho')$
are those where $V_\lambda\ra V$, $\lambda\in L$ is a covering family of schemes with respect to the fppf topology.
The \emph{structure sheaf} $\O_{\shM}$ is given by the formula
\begin{equation}
\label{structure sheaf stack}
\Gamma(\hat{V},\O_{\shM})=\Gamma(V',\O_{V'})=\Gamma(V,\O_V) 
\end{equation} 
for $\hat{V}=(V,\shE',\rho')$, with  obvious notion of restriction maps. 

To define the \emph{tautological sheaf} $\shT_\shM$ on the stack, first note that for each   $\hat{V}=(V,\shE',\rho')$, the 
splitting gerbe $V'$ comes with two locally free sheaves, 
namely $\shE'$ and the tautological sheaf $\shF_{V'}^\taut$, both of rank $n=\deg(\Lambda^\azu)$ and weight one.
In turn, the  Hom sheaf $\uHom(\shF_{V'}^\taut,\shE')$ has weight zero, and thus can be seen as a locally free sheaf
$\shT_{\shM,\hat{V}}$ on $V$ of rank $n^2$.
This said, we define 
\begin{equation}
\label{tautological sheaf stack}
\Gamma(  \hat{V}, \shT_\shM) = \Hom(\shF_{V'}^\taut,\shE')= \Gamma(V, \shT_{\shM,\hat{V}}),
\end{equation} 
with   obvious notion of restriction maps. 
The rings \eqref{structure sheaf stack} act  by scalar multiplication, which turns $\shT_\shM$ into a presheaf of $\O_\shM$-modules.

\begin{proposition}
\mylabel{tautological sheaf rank weight}
The presheaf $\shT_\shM$  on the Artin stack $\shM$ satisfies the sheaf axiom. As $\O_\shM$-module, it is locally free of rank $n^2$
and pure of weight one.   
\end{proposition}

\proof
First note that $\GG_m$ acts on the objects $\hat{V}=(V,\shE',\psi')$ via scalar multiplication on the sheaf $\shE'$.
The induced action on \eqref{tautological sheaf stack} is again via scalar multiplication, which has weight $w=1$.
The formation of the Hom sheaves $\uHom(\shF_{V'}^\taut,\shE')$ 
is compatible with morphisms  $\hat{V}_1\ra \hat{V_2}$ in the category $\shM$, 
and this  implies the sheaf axiom for   $\shT_\shM$. 
Since $\uHom(\shF_{V'}^\taut,\shE')$ are locally free 
of rank $n^2$, the same holds for $\shT_\shM$.
\qed

\medskip
We call $\shT_\shM$ the \emph{tautological sheaf} on the Artin stack $\shM$. The crucial observation now is that
the local sections \eqref{tautological sheaf stack} also  carry the structure of an $\Lambda^\azu$-module:
For each object $\tilde{U}=(U,h,\shF,\varphi)$ from the splitting gerbe $V'$, we get
$$
\Lambda^\azu\otimes\O_U\stackrel{\varphi^{-1}}{\lra}\uEnd(\shF) \stackrel{\circ}{\lra} \uEnd(\shT_\shM|U').
$$
The map on the right arises from $\shF=\shF^\taut_{V'}|U$, and is given by 
left composition with  respect to the Hom sheaf $\shT_\shM|U'=\uHom(\shF^\taut_{V'}|U,\shE'|U)$.
The above   is compatible with  restrictions, and   turns $\shT_\shM$
into a sheaf  of $\Lambda^\azu\otimes_R\O_\shM$-modules. The latter can be seen as inclusion
$\Lambda^\azu\otimes_R\O_\shM\subset\uEnd_{\O_\shM}(\shT_\shM)$
of Azumaya  algebras over $\shM$, both of weight zero. We now set
$$
\shA_\shM=\uEnd_{\O_\shM}(\shT_\shM)\quadand \shA^0_\shM= \uEnd_{\Lambda^\azu\otimes\O_\shM}(\shT_\shM).
$$
In other words, $\shA^0_\shM\subset\shA_\shM$ is the \emph{commutant} of $\Lambda^\azu\otimes\O_\shM\subset\shA_\shM$.
Taking the \emph{bi-commutant} $\uEnd_{\shA^0_\shM}(\shT_\shM)$, we arrive at a commutative diagram
\begin{equation}
\label{three azumaya algebras}
\begin{tikzcd} 
\Lambda^\azu\otimes_R {\O_\shM} \ar[r]\ar[d]	& \shA_\shM\ar[d,"\id"]		& \shA^0_\shM\ar[l]\ar[d,"\id"]\\
\uEnd_{\shA^0_\shM}(\shT_\shM)\ar[r]	& \uEnd_{\O_\shM}(\shT_\shM)	& \uEnd_{\Lambda^\azu\otimes\O_\shM}(\shT_\shM)\ar[l],
\end{tikzcd}
\end{equation} 
where all maps are injective. The following   complements Theorem \ref{comparison functor equivalence}, 
and will elucidate the significance  of the tautological sheaf  $\shT_\shM$: 

\begin{lemma}
\mylabel{azumaya of weight zero}
In the above diagram, all terms are Azumaya algebras of weight zero, and the vertical arrow on the left is bijective.  
\end{lemma}

\proof
Obviously, $\Lambda^\azu\otimes_R {\O_\shM}$ and $\shA_\shM=\uEnd_{\O_\shM}(\shT_\shM)$ are  Azumaya algebras over $\O_\shM$.
It then follows from \cite{Auslander; Goldman 1960}, Theorem 3.3 that this carries over to the  commutant $ \shA^0_\shM$, and that the inclusion
$\Lambda^\azu\otimes_R {\O_\shM} \subset  \shA^0_\shM$ is an equality. Since $\shA_\shM$ has weight zero, the same holds for the quasicoherent subsheaves.
\qed

\medskip
Since the terms  have weight zero, the upper horizontal row in \eqref{three azumaya algebras}  
descends along the $\GG_m$-gerbe $ \shM \ra(\Aff/Q)$  to sheaves of locally free algebras
$$
\Lambda^\azu\otimes_R\O_Q\lra \shA_Q\longleftarrow \shA^0_Q
$$ 
on the algebraic space $Q$.  Now recall that an associative $R$-algebra $A$ is an Azumaya algebra
if and only if the underlying $R$-module is locally free of finite rank,
and  the canonical map $A\otimes_RA^\op\ra\End_R(A)$ given by left-right multiplication 
is bijective. Using this characterization, we infer that the above are Azumaya algebras over $Q$,
and give rise to  classes $[\Lambda^\azu\otimes_R\O_Q]$ and  $[\shA_Q]$ and $[\shA^0_Q]$ in the Brauer group
$\Br(Q)$.
By construction, $\shA_\shM=\shA_Q\otimes\O_\shM$ becomes the endomorphism algebra for the tautological sheaf   $\shT_\shM$,
and thus has trivial class in $\Br(\shM)$. Note that this observation does not carry over to $\shA_Q$, because the tautological sheaf
 $\shT_\shM$ has weight one and therefore does not descend.

Now recall that the quotient map $X\ra X/H^\op=Q$ is a torsor with respect to  $H_Q$ for  the twisted form $H=\Aut_{\Lambda^\azu/R}$
of $\PGL_n$, and that  $G=U_{\Lambda^\azu/R}$ is a twisted form of $\GL_n$. We have  a  short exact sequence 
$1\ra \GG_{m,Q}\ra G_Q\ra H_Q\ra 1$, and  the torsor class $[X]\in H^1(Q,H_Q)$ yields
via  the non-abelian coboundary some $\partial[X]\in H^2(Q,\GG_m)$.
As a direct consequence of the preceding result, we obtain:

\begin{theorem}
\mylabel{tautological classes}
In the  Brauer group  $\Br(Q)\subset H^2(Q,\GG_m)$ of the  algebraic space $Q=X/H^\op$, we have 
$$
[\shA_Q]=\partial [X]  \quadand [\shA^0_Q]=[\Lambda^\azu\otimes_R\O_Q].
$$
\end{theorem}

\proof
By definition of the non-abelian coboundary, the class $\partial[X]$ is represented by the $\GG_m$-gerbe  $[X/G_Q/Q]$.
In light of Theorem \ref{comparison functor equivalence}, the $\GG_m$-gerbe $\shM=\shM^{\Lambda^\azu}_{\Lambda/R}$ is another representative.
The tautological sheaf $\shT_\shM$ is locally free, of rank $n^2$ and weight one, and by construction
$\shA_Q\otimes\O_\shM=\uEnd(\shT_\shM)$.  Now de Jong's observation (\cite{de Jong 2006}, Lemma 2.14) gives $[\shA_Q]=[\shM]$.

For the second assertion, we first work   on the stack $\shM$. Set $T=\shT_\shM$.
The sheaf of Azumaya algebras $A=\End(T)$  contains   
both $B=\Lambda^\azu\otimes_R\O_\shM$ and $C=\shA^0_\shM$, and thus becomes a sheaf of  modules over $B\otimes_AC^\op$.
According to   Lemma \ref{azumaya of weight zero}, the canonical maps $B\ra\End_C(T)$ and $C\ra \End_B(T)$ are bijective.
Since these sheaves  have weight zero, the very same statements  hold  on the algebraic space $Q$: 
Changing notation, we set 
$$
T=\shT_Q\quadand A=\End(T)\quadand B=\Lambda^\azu\otimes_R\O_Q\quadand C=\shA^0_Q.
$$
Again the canonical maps $B\ra\End_C(T)$ and $C\ra \End_B(T)$ are bijective. The equality $[B]=[C]$ in the Brauer group $\Br(Q)$
follows (see for example \cite{Curtis; Reiner 1990}, Theorem 3.54).
\qed

%===========================================================

\end{document}